\RequirePackage{fix-cm}
 \documentclass[smallextended]{svjour3}
 \smartqed 

\usepackage{amssymb,amsmath,amsfonts,latexsym}
\usepackage{hyperref}

\usepackage{mathtools}
\usepackage{multirow,array}
\usepackage{graphicx}
\usepackage{subfigure}
\usepackage{epstopdf}
\usepackage{float}
\usepackage{color, xcolor}
\usepackage{comment}
\usepackage{url}
\usepackage[subrefformat=parens]{subcaption}
\usepackage{mathtools}
\usepackage[section]{placeins}
\usepackage{amsmath}
\usepackage{bm}

\allowdisplaybreaks[3]

\def\diam{\mathop{\mathrm{diam}}}
\def\diag{\mathop{\mathrm{diag}}}

\def\div{\mathop{\mathrm{div}}}
\def\divh{\mathop{\mathrm{div}_h}}

\def\<{\mathop{\textless}}
\def\>{\mathop{\textgreater}}

\def\star{\mathop{\rm{star}}}

\spnewtheorem{thr}{Theorem}{\bf}{\it}
\spnewtheorem{defi}{Definition}{\bf}{\it}
\spnewtheorem{lem}{Lemma}{\bf}{\it}
\spnewtheorem{coro}{Corollary}{\bf}{\it}
\spnewtheorem{assume}{Assumption}{\bf}{\it}
\spnewtheorem{ex}{Example}{\bf}{\it}
\spnewtheorem{Case}{Case}{\bf}{\it}
\spnewtheorem*{pf*}{Proof}{\bf}{\rm}
\spnewtheorem*{rem*}{Remark:}{\it}{\it}
\spnewtheorem*{ex*}{Example:}{\it}{\it}
\spnewtheorem{Cond}{Condition}{\bf}{\it}
\spnewtheorem{rem}{Remark}{\it}{\it}

\spnewtheorem*{lem1*}{Lemma 1}{\bf}{\rm}
\spnewtheorem*{appendix1*}{Appendix A}{\bf}{\rm}

\newcounter{sone}
\newcounter{stwo}
\newcounter{sthree}
\newcounter{sfour}
\newcounter{sfive}
\newcounter{ssix}
\newcounter{lone}
\newcounter{ltwo}
\newcounter{lthree}
\newcounter{lfour}
\newcounter{lfive}
\newcounter{lsix}
\makeatletter
    
    \@addtoreset{equation}{section}

\journalname{}

\begin{document}

\title{Nitsche's method under a semi-regular mesh condition}

\titlerunning{Nitsche's method under a semi-regular mesh condition}

\author{
      Hiroki Ishizaka 
 }

\authorrunning{H. Ishizaka}

\institute{Hiroki Ishizaka
\at
Team FEM, Matsuyama, Japan \\
\email{h.ishizaka005@gmail.com}
}

\date{Received: date / Accepted: date}
\maketitle

\begin{abstract}
Nitsche's method is a numerical approach that weakly enforces boundary conditions for partial differential equations. In recent years, Nitsche's method has experienced a revival owing to its natural application in modern computational methods, such as the cut and immersed finite element methods. This study investigates Nitsche's methods based on an anisotropic weakly over-penalised symmetric interior penalty method for Poisson and Stokes equations on convex domains. As our primary contribution, we provide a new proof for the consistency term, which allows us to obtain an estimate of the anisotropic consistency error. The key idea of the proof is to apply the relationship between the Crouzeix and Raviart finite element space and the Raviart--Thomas finite element space.  We present the error estimates in the energy norm on anisotropic meshes. We compared the calculation results for the anisotropic mesh partitions in the numerical experiments.

\keywords{Nitsche's method \and Poisson and Stokes problems \and CR finite element method \and RT finite element method \and Semi-regular mesh condition}

\subclass{65D05 \and 65N30}

\end{abstract} 

\section{Introduction} \label{sec:1}

\setcounter{section}{1} \setcounter{equation}{0}
In this study, we examine Nitsche's methods for the Poisson and Stokes equations with a non-homogeneous Dirichlet boundary condition. The Nitsche method, introduced by Joachim A. Nitsche in 1971 (\cite{Nt71}), is a numerical technique for the weak imposition of Dirichlet boundary conditions, particularly in finite element analysis. This method is a penalty method that adds appropriate boundary terms to the variational formulation. This formulation can be interpreted as a special case of classical interior penalty discontinuous Galerkin (dG) methods (e.g., \cite{ErnGue21b,PieErn12}).
Nitsche's method is extensively utilised in fluid-structure interaction problems, contact problems, and in conjunction with the discontinuous Galerkin method. Numerous references exist concerning the theoretical and empirical analyses of Nitsche's method, as well as a priori and a posteriori error analysis, for instance, \cite{ChoErnPig20,ChoHil13,ChoHil13b,ChoHilLleRen22,DabDel21,DabMarVoh20a,DabMarVoh20b}. However, the analysis traditionally imposes the shape-regularity mesh condition, which is standard in finite element error analysis. Under the shape-regular family of triangulations, triangles or tetrahedra must not become excessively flat. Conversely, anisotropic meshes, where elements are not uniform in all directions, prove effective in problems where the solution demonstrates anisotropic behaviour in certain domain directions. Nonetheless, developing accurate and efficient finite-element schemes to solve partial differential equations across various domains remains challenging. In \cite{Ish21,Ish24b,IshKobTsu21a,IshKobTsu21c}, we introduced a novel parameter (Definition \ref{defi1}) and a geometric condition (Assumption \ref{neogeo=assume}) to supplant the maximum angle condition. This new condition can be readily computed numerically using the finite element method, potentially benefiting a posteriori error analyses. In this study, we conduct an anisotropic error analysis of Nitsche's method employing this new condition.

\color{black}
Throughout this study, let $\Omega \subset \mathbb{R}^d$, $d \in \{ 2 , 3 \}$ be a bounded polyhedral domain. Furthermore, we assume that $\Omega$ is convex if necessary. Let $f \in L^2(\Omega)$ and $g \in H^{\frac{1}{2}}(\partial \Omega)$. We consider the Poisson problem with a nonhomogeneous Dirichlet boundary condition. Find $u: \Omega \to \mathbb{R}$ such that
\begin{align}
\displaystyle
- \varDelta u  = f \quad \text{in $\Omega$}, \quad u = g \quad \text{on $\partial \Omega$}.\label{ell_eq}
\end{align}
The variational formulation of the Poisson problem \eqref{ell_eq} is to find $u \in V(g)$ such that
\begin{align}
\displaystyle
a(u,v) \coloneq ({\bm \nabla} u , {\bm \nabla} v)_{\Omega} = (f , v)_{\Omega} \quad \forall v \in V(0), \label{ell_weak}
\end{align}
where $(\cdot,\cdot)_{\Omega}$ denotes the $L^2$-scalar product over $\Omega \subset \mathbb{R}^d$, 
\begin{align*}
\displaystyle
V(g) \coloneq \left \{ v \in H^1(\Omega): \ v = g \text{ on } \partial \Omega \right \}.
\end{align*}
The variational problem \eqref{ell_weak} is well-posed, e.g., see \cite[Proposition 31.12]{ErnGue21b}. {We consider the classical Lagrange finite element method and a symmetric penalty method for this problem.}

Let $\mathbb{T}_h = \{ T \}$ be a simplicial mesh of $\overline{\Omega}$ composed of closed $d$-simplices, such as $\overline{\Omega} = \bigcup_{T \in \mathbb{T}_h} T$, where $h \coloneq \max_{T \in \mathbb{T}_h} h_{T}$ and $ h_{T} \coloneq \diam(T)$. For simplicity, we assume that $\mathbb{T}_h$ is conformal, that is, $\mathbb{T}_h$ is a simplicial mesh of $\overline{\Omega}$ without hanging nodes. For $k \in \mathbb{N}_0 \coloneq \mathbb{N} \cup \{ 0 \}$, $\mathbb{P}^k(T)$ is spanned by the restriction of $T$ by polynomials in $\mathbb{P}^k$, where  $\mathbb{P}^k$ denotes the space of polynomials with a maximum of $k$ degrees. Let $P_{dc,h}^{1}$ be a discontinuous $\mathbb{P}^1$ space, defined in \eqref{dis=sp}. For Nitsche's method, we define a finite element space as $P_{c,h}^{1} \coloneq \{ v_h \in P_{dc,h}^{1}: \ [\![ v_h ]\!]_F = 0 \ \forall F \in \mathcal{F}_h^i \}$, where  $[\![ v_h ]\!]_F$ is the jump of $v_h$ across $F$, and $\mathcal{F}_h^i$ is the set of interior faces. The formulation is obtained by removing the consistency terms at all interfaces in classical dG schemes. The symmetric Nitsche formulation with $g \equiv 0$ obtained is as follows: Find $u_h \in P_{c,h}^{1}$ such that
\begin{align}
\displaystyle
a_h^{sip}(u_h , v_h) = (f,v_h)_{\Omega} \quad \forall {v_h} \in P_{c,h}^{1}, \label{Sym=Nit1}
\end{align}
where $P_{c,h,*}^1 \coloneq  H^2(\Omega) + P_{c,h}^1$, and $a_h^{sip}: P_{c,h,*}^{1} \times P_{c,h}^{1} \to \mathbb{R}$ is defined as 
\begin{align*}
\displaystyle
a_h^{sip}(v,w_h) &\coloneq \int_{\Omega} {\bm \nabla} v \cdot {\bm \nabla} w_h dx + \sum_{F \in \mathcal{F}_h^{\partial}} \gamma^{sip} \kappa_{F(0)}  \int_F v w_hds \\
&\quad - \sum_{F \in \mathcal{F}_h^{\partial}} \int_F \left( \bm \nabla v \cdot \bm n  w_h  + v  \bm \nabla w_h  \cdot \bm n \right ) ds,
\end{align*}
where $\mathcal{F}_h^{\partial}$ denotes the set of faces on the boundary $\partial \Omega$, $\kappa_{F(0)}$ is defined in Section \ref{sec=notation} and the penalty parameter $\gamma^{sip} \> 0$ is chosen so that the coercivity is established. We define a finite-element space as $P_{c,h}^L := \{ v_h \in P_{c,h}^1; \ v_h |_{\partial \Omega} = 0 \}$. The Lagrange formulation with $g \equiv 0$ obtained is as follows: Find $u_h^L \in P_{c,h}^{L}$ such that
\begin{align}
\displaystyle
\left( a_h^{sip}(u_h^L , v_h^L) = \right) \int_{\Omega} \bm \nabla u_h^L \cdot \bm \nabla v_h^L dx = (f,v_h^L)_{\Omega} \quad \forall v_h^L \in P_{c,h}^{L}. \label{Lag1}
\end{align}
Let $d=2$ and let $I_T^L: H^2(T) \to \mathbb{P}^1(T)$ be the typical Lagrange interpolation operator. Reference \cite[Corollary 1]{IshKobTsu21c} presented the following anisotropic interpolation error under the semi-regular condition \eqref{NewGeo}: 
\begin{align}
\displaystyle 
|u - I_T^L u|_{H^1(T)} &\leq c \sum_{i=1}^2 h_i \left | \frac{\partial u}{\partial \bm {r_i}} \right |_{H^{1}(T)}, \label{Lag=int}
\end{align}
where $h_i$ and $\bm r_i$ are as described in Section \ref{sec:4}. We use the above anisotropic interpolation error and the Galerkin orthogonality to derive the error estimate of the Lagrange scheme \eqref{Lag1}. Let $I_h^L: H^2(\Omega) \to P_{c,h}^1$ be the global Lagrange interpolation operator. For the solution $u \in H^2(\Omega) \cap H_0^1(\Omega)$ of the problem \eqref{ell_weak}, we have
\begin{align}
\displaystyle
|u - u_h^L|_{H^1(\Omega)}
&\leq c \inf_{w_h \in  P_{c,h}^{L}} |u - w_h|_{H^1(\Omega)} \leq c  |u - I_h^L u|_{H^1(\Omega)} \notag \\
&\leq c \sum_{T \in \mathbb{T}_h} \sum_{i=1}^2 h_i \left | \frac{\partial u}{\partial \bm {r_i}} \right |_{H^{1}(T)}. \label{Lag2}
\end{align}
This error estimate may be optimal on anisotropic meshes. 

On the other hand, we establish consistency, boundedness and discrete coercivity for the bilinear form to obtain a convergence analysis for the symmetric Nitsche formulation \eqref{Sym=Nit1}. We define a norm as $ | v |_{SN} := ( | v |^2_{H^1(\Omega)} + \sum_{F \in \mathcal{F}_h^{\partial}} \kappa_{F(0)} \| v \|_{L^2(F)}^2 )^{1/2}$ for any $v \in P_{c,h,*}^1$. By choosing the parameter $\gamma^{sip} \> 0$ appropriately, the following coersivity holds: there eixists $\alpha \> 0$ such that
\begin{align*}
\displaystyle
a_h^{sip}(w_h,w_h) \geq \alpha | w_h |_{SN}^2 \quad \forall w_h \in  P_{c,h}^1.
\end{align*}
Thus, we have the discrete stability:
\begin{align}
\displaystyle
| w_h |_{SN} \leq \frac{1}{\alpha} \sup_{v_h \in P_{c,h}^1} \frac{a_h^{sip}(w_h,v_h)}{| v_h |_{SN}}. \label{sta=SymNit}
\end{align}
By the integration by parts,
\begin{align*}
\displaystyle
a_h^{sip}(u, v_h)
&= \int_{\Omega} \bm \nabla u \cdot  \bm \nabla v_h dx - \sum_{F \in \mathcal{F}_h^{\partial}} \int_F \bm \nabla u \cdot \bm n  v_h ds = \int_{\Omega} - \varDelta u  v_h dx \\
&= \int_{\Omega} f v_h dx = a_h^{sip}(u_h , v_h),
\end{align*}
which leads to
\begin{align}
\displaystyle
a_h^{sip}(u - u_h, v_h) = 0 \quad \forall v_h \in P_{c,h}^1. \label{Galerkin=ortho}
\end{align}
Therefore, using the discrete stability and \eqref{Galerkin=ortho}, 
\begin{align}
\displaystyle
|u - u_h|_{SN}
&\leq |u - I_h^L u|_{SN} + |I_h^L u - u_h|_{SN} \notag\\
&\leq |u - I_h^L u|_{SN} +  \frac{1}{\alpha} \sup_{v_h \in P_{c,h}^1} \frac{ a_h^{sip}(I_h^L u - u,v_h)}{| v_h |_{SN}}. \label{error=SymNit}
\end{align}
We estimate the consistency term.  For $F \in \mathcal{F}_h^{\partial}$ with $F = \partial T \cap \partial \Omega$ and any $(v,v_h) \in P_{c,h,*}^{1} \times P_{c,h}^{1}$, the trace (\cite[Lemma 1]{Ish24}) and the H\"older inequalities yield
\begin{align*}
\displaystyle
\left| \int_F \bm \nabla v \cdot \bm n v_h  ds \right|
&\leq c \left(  h \| \bm \nabla v \|_{L^2(T)^d} + h^{\frac{3}{2}} \| \bm \nabla v \|_{L^2(T)^d}^{\frac{1}{2}} | \bm \nabla v |_{H^1(T)^d}^{\frac{1}{2}} \right) \\
&\quad \times \kappa_{F(0)}^{\frac{1}{2}} \|  v_{h} \|_{L^2(F)}.
\end{align*}
We set $v := u - I_T^L u$. The interpolation error \eqref{Lag=int} leads to
\begin{align} 
\displaystyle  
&h \| \bm \nabla (u - I_T^L u) \|_{{L^2(T)^2}} + h^{\frac{3}{2}}  \| \bm \nabla (u - I_T^L u) \|_{{L^2(T)^2}}^{\frac{1}{2}} | \bm \nabla (u - I_T^L u)|_{{H^1(T)^2}}^{\frac{1}{2}} \notag  \\  
&\quad \leq c h {\sum_{i=1}^2} h_i \left | \frac{\partial u}{\partial \bm r_i} \right |_{H^{1}(T)} + c h^{\frac{3}{2}} \left( {\sum_{i=1}^2} h_i \left | \frac{\partial u}{\partial \bm r_i} \right |_{H^{1}(T)} \right)^{\frac{1}{2}} |u|_{H^2(T)}^{\frac{1}{2}} \label{SN1} \\
&\quad \leq c h^{\frac{1}{2}} \left( {\sum_{i=1}^2} h_i \left | \frac{\partial u}{\partial \bm r_i} \right |_{H^{1}(T)} \right)^{\frac{1}{2}} |u|_{H^2(T)}^{\frac{1}{2}} \quad \text{if $h \leq 1$}. \notag
\end{align}
Compared to the error estimate \eqref{Lag2} of the Lagrange scheme, the estimate of the consistency term in the symmetric Nitsche formulation may degenerate. The second term in \eqref{SN1} absorbs the first term on the right-hand side. Therefore, the second term may reduce the computational efficiency of Nische's method for anisotropic meshes.

In Section \ref{sec:3}, we introduce Nische's method \eqref{ell_dis} for addressing the Poisson equation, utilising a weakly over-penalised symmetric interior penalty (WOPSIP) approach. A notable characteristic of the WOPSIP method is the absence of consistency terms within the scheme, thereby circumventing the degeneracy of the error estimate previously described. \color{black}The WOPSIP method is a simple dG method similar to the classical Crouzeix--Raviart (CR) finite element method. Brenner et al. first proposed the WOPSIP method \cite{BarBre14,BreOweSun08}, and several further studies have considered similar techniques for anisotropic meshes \cite{Ish24,Ish24c}. Our formulation is non-conforming. Therefore, the error between the exact and finite-element approximation solutions for an energy norm can be divided into two parts. One part is the optimal approximation error in finite element spaces, and the other is the consistency error term. For the former, the classical CR interpolation error (Theorem \ref{DGCE=thr3}) is used. For estimating the consistency error term on anisotropic meshes, we use the relation between the CR finite element interpolation and the lowest order Raviart--Thomas (RT) finite element interpolation (Lemma \ref{Nit=lem1}). {This relation derives a consistency error (Lemma \ref{asy=con}).}

Furthermore, in Section \ref{sec:7}, we introduce the Nitsche method for the pressure-robust formulation of the Stokes equation on anisotropic meshes. In \cite{Lin14}, as {for} the Stokes element, the CR finite element space for {the} velocity and the elementwise constant discontinuous space for {the} pressure were proposed for the analysis. We extend the approach.

The remainder of this paper is organised as follows. In Section 2, we introduce {some notations} and a new scheme for the Poisson equation. Section 3 provides a new parameter, geometric condition, and anisotropic interpolation error estimates.  In Section 4, we present the discrete Poincar\'e inequality. In Section 5, we discuss stability and error estimates. Section 6 presents the Nitsche scheme for the Stokes equation and error estimates. Section 7 presents the numerical results.

Throughout this paper, we denote by $c$ a constant independent of $h$ (defined later) and the angles and aspect ratios of simplices, unless specified otherwise {all constants $c$ are bounded if the maximum angle is bounded}. These values vary across different contexts. The notation $\mathbb{R}_+$ denotes a set of positive real numbers.

\color{black}
\section{Nitsche method for the Poisson equation} \label{sec:3}
\setcounter{section}{2} \setcounter{equation}{0}
\subsection{{Notations}} \label{sec=notation}
\textbf{Mesh faces and jumps.} 
We set $\mathcal{F}_h := \mathcal{F}_h^i \cup \mathcal{F}_h^{\partial}$. For any $F \in \mathcal{F}_h$, we define the unit normal {{$\bm {n}_F$}} to $F$ as follows: (\roman{sone}) If  $F \in \mathcal{F}_h^i$ with $F = T_1 \cap T_2$, $T_1,T_2 \in \mathbb{T}_h$, let ${\bm {n}_1}$ and ${\bm {n}_2}$ be the outwards unit normals of $T_1$ and $T_2$, respectively. Then, ${\bm {n}_F}$ is either of $\{ {\bm {n}_1} , {\bm {n}_2}\}$; (\roman{stwo}) If $F \in \mathcal{F}_h^{\partial}$, ${\bm {n}_F}$ is the unit outwards normal ${\bm {n}}$ to $\partial \Omega$. We define the broken (piecewise) Sobolev space as
\begin{align*}
\displaystyle
H^1(\mathbb{T}_h) &\coloneq \left\{ \varphi \in L^2(\Omega); \ \varphi|_{T} \in H^1(T) \ \forall T \in \mathbb{T}_h  \right\}
\end{align*}
with the seminorm
\begin{align*}
\displaystyle
| \varphi |_{H^1(\mathbb{T}_h)} &\coloneq \left( \sum_{T \in \mathbb{T}_h} | \varphi |^2_{H^1(T)} \right)^{\frac{1}{2}} \quad \varphi \in H^1(\mathbb{T}_h).
\end{align*}
Let $\varphi \in H^1(\mathbb{T}_h)$. Suppose that $F \in \mathcal{F}_h^i$ with $F = T_1 \cap T_2$, $T_1,T_2 \in \mathbb{T}_h$. We set $\varphi_1 \coloneq \varphi{|_{T_1}}$ and $\varphi_2 \coloneq \varphi{|_{T_2}}$. We set two nonnegative real numbers $\omega_{T_1,F}$ and $\omega_{T_2,F}$ such that
\begin{align*}
\displaystyle
\omega_{T_1,F} + \omega_{T_2,F} = 1.
\end{align*}
The jump and skew-weighted averages of $\varphi$ across $F$ are then defined as
\begin{align*}
\displaystyle
[\![\varphi]\!] \coloneq [\! [ \varphi ]\!]_F := \varphi_1 - \varphi_2, \quad  \{\! \{ \varphi\} \! \}_{\overline{\omega}} \coloneq  \{\! \{ \varphi\} \! \}_{\overline{\omega},F} \coloneq \omega_{T_2,F} \varphi_1 + \omega_{T_1,F} \varphi_2.
\end{align*}
For a boundary face $F \in \mathcal{F}_h^{\partial}$ with $F = \partial T \cap \partial \Omega$, $[\![\varphi ]\!]_F \coloneq \varphi|_{T}$ and $\{\! \{ \varphi \} \!\}_{\overline{\omega}} \coloneq \varphi |_{T}$. For any ${\bm {v}} \in H^{1}(\mathbb{T}_h)^d$, the notation
\begin{align*}
\displaystyle
&[\![{\bm {v}} \cdot {\bm {n}}]\!] \coloneq [\![ {\bm {v}} \cdot {\bm {n}} ]\!]_F \coloneq {\bm {v}_1} \cdot {\bm {n}_F} - {\bm {v}_2} \cdot {\bm {n}_F},\\
&  \{\! \{ {\bm {v}}\} \! \}_{\omega} \coloneq  \{\! \{ {\bm {v}} \} \! \}_{\omega,F} \coloneq \omega_{T_1,F} {\bm {v}_1} + \omega_{T_2,F} {\bm {v}_2}
\end{align*}
for the jump in the normal component and the weighted average of ${\bm {v}}$. For any ${\bm {v}} \in H^{1}(\mathbb{T}_h)^d$ and $\varphi \in H^{1}(\mathbb{T}_h)$, we have that
\begin{align*}
\displaystyle
[\![ ({\bm {v}} \varphi) \cdot {\bm {n}} ]\!]_F
&=  \{\! \{ {\bm {v}} \} \! \}_{\omega,F} \cdot {\bm {n}_F} [\! [ \varphi ]\!]_F + [\![ {\bm {v}} \cdot {\bm {n}} ]\!]_F \{\! \{ \varphi\} \! \}_{\overline{\omega},F}.
\end{align*}

For any $F \in \mathcal{F}_h^{\partial}$, we define the $L^2$-projection $\Pi_F^{0}: L^2(F) \to \mathbb{P}^{0}(F)$ as
\begin{align*}
\displaystyle
\int_F (\Pi_F^{0} \varphi - \varphi)   ds = 0 \quad \forall \varphi \in L^2(F).
\end{align*}

For $\varphi \in H^1(\mathbb{T}_h)$, the broken gradient ${\bm\nabla_h}:H^1(\mathbb{T}_h) \to L^2(\Omega)^{d}$ is defined as
\begin{align*}
\displaystyle
({\bm\nabla_h} \varphi)|_T &\coloneq {\bm\nabla} (\varphi|_T) \quad \forall T \in \mathbb{T}_h.
\end{align*}

\textbf{Finite element spaces.} For $s \in \mathbb{N}_0$, we define the standard discontinuous finite-element space as
\begin{align}
\displaystyle
P_{dc,h}^{s} &\coloneq \left\{ p_h \in {L^{\infty}(\Omega)}; \ p_h|_{T} \in \mathbb{P}^{s}({T}) \quad \forall T \in \mathbb{T}_h  \right\}. \label{dis=sp}
\end{align}
The CR finite element space is defined as
\begin{align}
\displaystyle
V_{h}^{CR} &\coloneq  \biggl \{ \varphi_h \in P_{dc,h}^1: \  \int_F [\![ \varphi_h ]\!] ds = 0 \ \forall F \in \mathcal{F}_h^i \biggr \}. \label{CRdef}
\end{align}
Furthermore, we set
\begin{align*}
\displaystyle
| v_h |_{1,\beta} \coloneq \left( |v_h|_{H^1(\mathbb{T}_h)}^2 + |v_h|_{\beta}^2  \right)^{\frac{1}{2}}, \quad  |v_h|_{\beta} \coloneq \left( \sum_{F \in \mathcal{F}_h^{\partial}} \kappa_{F(\beta)} \| \Pi_F^0 v_h \|^2_{L^2(F)} \right)^{\frac{1}{2}}
\end{align*}
for any $v_h \in V_h^{CR}$ and $\beta \in \mathbb{R}_+$. Here, the  penalty parameter $\kappa_{F(\beta)}$ is defined as
\begin{align*}
\displaystyle
\kappa_{F(\beta)} \coloneq h^{- 2 \beta} \ell_{T,F}^{- 1}, \quad \ell_{T,F} \coloneq \frac{d! |T|_d}{|F|_{d-1}} \quad F \in \mathcal{F}_h^{\partial} \text{ with } F = \partial T \cap \partial \Omega, \ T \in \mathbb{T}_h,
\end{align*}
{see \cite[Section 2.4]{Ish24}.} Here, $|\cdot|_d$ denotes the $d$-dimensional Hausdorff measure. For any $v_h \in V_h^{CR}$, $ | v_h |_{0} \leq  | v_h |_{1}$ for $h \leq 1$. 

\subsection{New scheme} \label{new=sch}
We introduce a new scheme for the Poisson equation as follows: For all $f \in L^2(\Omega)$ and $g \in H^{\frac{1}{2}}(\partial \Omega)$, we aim to determine $u_h \in V_h^{CR}$ such that
\begin{align}
\displaystyle
a_h(u_h , v_h) = \ell_h(v_h) \quad \forall v_h \in V_h^{CR}, \label{ell_dis}
\end{align}
where $a_h: (H^1(\Omega) + V_h^{CR}) \times (H^1(\Omega) + V_h^{CR})  \to \mathbb{R}$ and $\ell_h: H^1(\Omega) + V_h^{CR} \to \mathbb{R}$ are defined as
\begin{align*}
\displaystyle
a_h(u_h , v_h)
&\coloneq ( {\bm \nabla_h} u_h , {\bm \nabla_h} v_h)_{\Omega} + \sum_{F \in \mathcal{F}_h^{\partial}} \kappa_{F(1)} \langle \Pi_F^0 u_h  , \Pi_F^0 v_h \rangle_F,\\
\ell_h(v_h)
&\coloneq (f,v_h)_{\Omega} + \sum_{F \in \mathcal{F}_h^{\partial}} \kappa_{F(1)} \langle \Pi_F^0  g, \Pi_F^0 v_h \rangle_F,
\end{align*}
where $\langle \cdot,\cdot \rangle_{F}$ denotes the $L^2$-scalar product over $F$. 

Using H\"older's inequality, we obtain
\begin{align*}
\displaystyle
|a_h(u_h , v_h)|
&\leq c | u_h |_{1,1} | v_h |_{1,1} \quad \forall u_h,v_h \in H^1(\Omega) + V_h^{CR}.
\end{align*}
Furthermore, it holds that
\begin{align*}
\displaystyle
a_h(v_h , v_h) = | v_h |_{1,1}^2 \quad \forall v_h \in H^1(\Omega) + V_h^{CR}.
\end{align*}

\begin{rem}
Section \ref{sec:5} proves the discrete Poincar\'e inequality. The stability of the scheme and $\ell_h \in (V_h^{CR})^{\prime}$ are shown via its discrete inequality. From the Lax--Milgram lemma, the discrete problem \eqref{ell_dis} is well-posed.
\end{rem}

\color{black}
\section{Anisotropic interpolation error estimates}  \label{sec:4}

\setcounter{section}{3} \setcounter{equation}{0} 
Our strategy for interpolation errors on anisotropic meshes was proposed by \cite{Ish21,IshKobTsu21a,IshKobTsu21c}.

\subsection{Reference Elements} \label{reference}
We first define the reference elements $\widehat{T} \subset \mathbb{R}^d$.

\subsubsection*{Two-dimensional case} \label{reference2d}
Let $\widehat{T} \subset \mathbb{R}^2$ be a reference triangle with vertices $\hat{\bm p}_1 := (0,0){^{\top}}$, $\hat{\bm p}_2 := (1,0){^{\top}}$, and $\hat{\bm p}_3 := (0,1){^{\top}}$. 

\subsubsection*{Three-dimensional case} \label{reference3d}
In the three-dimensional case, we consider the following two cases: (\roman{sone}) and (\roman{stwo}); see Condition \ref{cond2}.

Let $\widehat{T}_1$ and $\widehat{T}_2$ be reference tetrahedra with the following vertices:
\begin{description}
   \item[(\roman{sone})] $\widehat{T}_1$ has vertices $\hat{\bm p}_1 := (0,0,0){^{\top}}$, $\hat{\bm p}_2 := (1,0,0){^{\top}}$, $\hat{\bm p}_3 := (0,1,0){^{\top}}$, and $\hat{\bm p}_4 := (0,0,1)^T$;
 \item[(\roman{stwo})] $\widehat{T}_2$ has vertices $\hat{\bm p}_1 := (0,0,0){^{\top}}$, $\hat{\bm p}_2 := (1,0,0){^{\top}}$, $\hat{\bm p}_3 := (1,1,0){^{\top}}$, and $\hat{\bm p}_4 := (0,0,1)^T$.
\end{description}
Therefore, we set $\widehat{T} \in \{ \widehat{T}_1 , \widehat{T}_2 \}$. 
Note that the case (\roman{sone}) is called \textit{the regular vertex property}, see \cite{AcoDur99}.

\subsection{Two-step Affine Mapping} \label{element=cond}
To an affine simplex $T \subset \mathbb{R}^d$, we construct two affine mappings $\Phi_{\widetilde{T}}: \widehat{T} \to \widetilde{T}$ and $\Phi_{T}: \widetilde{T} \to T$. First, we define the affine mapping $\Phi_{\widetilde{T}}: \widehat{T} \to \widetilde{T}$ as
\begin{align}
\displaystyle
\Phi_{\widetilde{T}}: \widehat{T} \ni \hat{\bm x} \mapsto \tilde{\bm x} := \Phi_{\widetilde{T}}(\hat{\bm x}) := {A}_{\widetilde{T}} \hat{\bm x} \in  \widetilde{T}, \label{aff=1}
\end{align}
where ${A}_{\widetilde{T}} \in \mathbb{R}^{d \times d}$ is an invertible matrix. We then define the affine mapping $\Phi_{T}: \widetilde{T} \to T$ as follows:
\begin{align}
\displaystyle
\Phi_{T}: \widetilde{T} \ni \tilde{\bm x} \mapsto \bm x := \Phi_{T}(\tilde{\bm x}) := {A}_{T} \tilde{\bm x} + \bm {b_{T}} \in T, \label{aff=2}
\end{align}
where $\bm {b_{T}} \in \mathbb{R}^d$ is a vector and ${A}_{T} \in O(d)$ denotes the rotation and mirror-imaging matrix. We define the affine mapping $\Phi: \widehat{T} \to T$ as
\begin{align*}
\displaystyle
\Phi := {\Phi}_{T} \circ {\Phi}_{\widetilde{T}}: \widehat{T} \ni \hat{\bm x} \mapsto \bm x := \Phi (\hat{\bm x}) =  ({\Phi}_{T} \circ {\Phi}_{\widetilde{T}})(\hat{\bm x}) = {A} \hat{\bm x} +\bm{ b_{T}} \in T, 
\end{align*}
where ${A} := {A}_{T} {A}_{\widetilde{T}} \in \mathbb{R}^{d \times d}$.

\subsubsection*{Construct mapping $\Phi_{\widetilde{T}}: \widehat{T} \to \widetilde{T}$} \label{sec221} 
We consider the affine mapping \eqref{aff=1}. We define the matrix $ {A}_{\widetilde{T}} \in \mathbb{R}^{d \times d}$ as follows. We first define the diagonal matrix as
\begin{align}
\displaystyle
\widehat{A} :=  \diag (h_1,\ldots,h_d), \quad h_i \in \mathbb{R}_+ \quad \forall i,\label{aff=3}
\end{align}
where $\mathbb{R}_+$ denotes the set of positive real numbers.

For $d=2$, we define the regular matrix $\widetilde{A} \in \mathbb{R}^{2 \times 2}$ as
\begin{align}
\displaystyle
\widetilde{A} :=
\begin{pmatrix}
1 & s \\
0 & t \\
\end{pmatrix}, \label{aff=4}
\end{align}
with the parameters
\begin{align*}
\displaystyle
s^2 + t^2 = 1, \quad t \> 0.
\end{align*}
For the reference element $\widehat{T}$, let $\mathfrak{T}^{(2)}$ be a family of triangles.
\begin{align*}
\displaystyle
\widetilde{T} &= \Phi_{\widetilde{T}}(\widehat{T}) = {A}_{\widetilde{T}} (\widehat{T}), \quad {A}_{\widetilde{T}} := \widetilde {A} \widehat{A}
\end{align*}
with the vertices $\tilde{\bm p}_1 := (0,0)^{\top}$, $\tilde{\bm p}_2 := (h_1,0)^{\top}$ and $\tilde{\bm p}_3 :=(h_2 s , h_2 t)^{\top}$. Then, $h_1 = |\tilde{\bm p}_1 - \tilde{\bm p}_2| \> 0$ and $h_2 = |\tilde{\bm p}_1 - \tilde{\bm p}_3| \> 0$. 

For $d=3$, we define the regular matrices $\widetilde{A}_1, \widetilde{A}_2 \in \mathbb{R}^{3 \times 3}$ as follows:
\begin{align}
\displaystyle
\widetilde{A}_1 :=
\begin{pmatrix}
1 & s_1 & s_{21} \\
0 & t_1  & s_{22}\\
0 & 0  & t_2\\
\end{pmatrix}, \
\widetilde{A}_2 :=
\begin{pmatrix}
1 & - s_1 & s_{21} \\
0 & t_1  & s_{22}\\
0 & 0  & t_2\\
\end{pmatrix} \label{aff=5}
\end{align}
with the parameters
\begin{align*}
\displaystyle
\begin{cases}
s_1^2 + t_1^2 = 1, \ s_1 \> 0, \ t_1 \> 0, \ h_2 s_1 \leq h_1 / 2, \\
s_{21}^2 + s_{22}^2 + t_2^2 = 1, \ t_2 \> 0, \ h_3 s_{21} \leq h_1 / 2.
\end{cases}
\end{align*}
Therefore, we set $\widetilde{A} \in \{ \widetilde{A}_1 , \widetilde{A}_2 \}$. For the reference elements $\widehat{T}_i$, $i=1,2$, let $\mathfrak{T}_i^{(3)}$, $i=1,2$, be a family of tetrahedra.
\begin{align*}
\displaystyle
\widetilde{T}_i &= \Phi_{\widetilde{T}_i} (\widehat{T}_i) =  {A}_{\widetilde{T}_i} (\widehat{T}_i), \quad {A}_{\widetilde{T}_i} := \widetilde {A}_i \widehat{A}, \quad i=1,2,
\end{align*}
with the vertices
\begin{align*}
\displaystyle
&\tilde{\bm p}_1 := (0,0,0)^{\top}, \ \tilde{\bm p}_2 := (h_1,0,0)^{\top}, \ \tilde{\bm p}_4 := (h_3 s_{21}, h_3 s_{22}, h_3 t_2)^{\top}, \\
&\begin{cases}
\tilde{\bm p}_3 := (h_2 s_1 , h_2 t_1 , 0)^{\top} \quad \text{for case (\roman{sone})}, \\
\tilde{\bm p}_3 := (h_1 - h_2 s_1, h_2 t_1,0)^{\top} \quad \text{for case (\roman{stwo})}.
\end{cases}
\end{align*}
Subsequently, $h_1 = |\tilde{\bm p}_1 - \tilde{\bm p}_2| \> 0$, $h_3 = |\tilde{\bm p}_1 - \tilde{\bm p}_4| \> 0$, and
\begin{align*}
\displaystyle
h_2 =
\begin{cases}
|\tilde{\bm p}_1 - \tilde{\bm p}_3| \> 0  \quad \text{for case (\roman{sone})}, \\
|\tilde{\bm p}_2 - \tilde{\bm p}_3| \> 0  \quad \text{for case (\roman{stwo})}.
\end{cases}
\end{align*}

\subsubsection*{Construct mapping $\Phi_{T}: \widetilde{T} \to T$}  \label{sec322}
We determine the affine mapping \eqref{aff=2} as follows. Let ${T} \in \mathbb{T}_h$ have vertices ${p}_i$ ($i=1,\ldots,d+1$). Let $b_{T} \in \mathbb{R}^d$ be the vector and ${A}_{T} \in O(d)$ be the rotation and mirror imaging matrix such that
\begin{align*}
\displaystyle
\bm p_{i} = \Phi_T (\tilde{\bm p}_i) = {A}_{T} \tilde{\bm p}_i + \bm{b_T}, \quad i \in \{1, \ldots,d+1 \},
\end{align*}
where vertices $\bm p_{i}$ ($i=1,\ldots,d+1$) satisfy the following conditions:

\color{black}
\begin{Cond}[Case in which $d=2$] \label{cond1}
Let ${T} \in \mathbb{T}_h$ have vertices ${\bm p_i}$ ($i=1,\ldots,3$). We assume that $\overline{{{\bm p}_2 {\bm p}_3}}$ is the longest edge of ${T}$, that is, $ h_{{T}} := |{{\bm p}_2 - {\bm p}_ 3}|$. We set $h_1 = |{{\bm p}_1 - {\bm p}_2}|$ and $h_2 = |{{\bm p}_1 - {\bm p}_3}|$. We then assume that $h_2 \leq h_1$. {Because $\frac{1}{2} h_T < h_1 \leq h_T$, ${h_1 \approx h_T}$.} 
\end{Cond}

\begin{Cond}[Case in which $d=3$] \label{cond2}
Let ${T} \in \mathbb{T}_h$ have vertices ${{\bm p}}_i$ ($i=1,\ldots,4$). Let ${L}_i$ ($1 \leq i \leq 6$) be the edges of ${T}$. We denote by ${L}_{\min}$  the edge of ${T}$ with the minimum length; that is, $|{L}_{\min}| = \min_{1 \leq i \leq 6} |{L}_i|$. We set $h_2 := |{L}_{\min}|$ and assume that 
\begin{align*}
\displaystyle
&\text{the endpoints of ${L}_{\min}$ are either $\{ {{\bm p}_1 , {\bm p}_3} \}$ or $\{ {{\bm p}_2 , {\bm p}_3}\}$}.
\end{align*}
Among the four edges sharing an endpoint with ${L}_{\min}$, we consider the longest edge ${L}^{({\min})}_{\max}$. Let ${{\bm p}}_1$ and ${\bm p}_2$ be the endpoints of edge ${L}^{({\min})}_{\max}$. Thus, we have
\begin{align*}
\displaystyle
h_1 = |{L}^{(\min)}_{\max}| = | {{\bm p}_1 - {\bm p}_2}|.
\end{align*}
We consider cutting $\mathbb{R}^3$ with a plane that contains the midpoint of the edge ${L}^{(\min)}_{\max}$ and is perpendicular to the vector ${{\bm p}_1 - {\bm p}_2}$. Thus, there are two cases. 
\begin{description}
  \item[(Type \roman{sone})] ${\bm p}_3$ and ${\bm p}_4$  belong to the same half-space;
  \item[(Type \roman{stwo})] ${\bm p}_3$ and ${\bm p}_4$  belong to different half-spaces.
\end{description}
In each case, we set
\begin{description}
  \item[(Type \roman{sone})] ${\bm p}_1$ and ${\bm p}_3$ as the endpoints of ${L}_{\min}$, that is, $h_2 =  |{\bm p}_1 - {\bm p}_3| $;
  \item[(Type \roman{stwo})] ${\bm p}_2$ and ${\bm p}_3$ as the endpoints of ${L}_{\min}$, that is, $h_2 =  |{\bm p}_2 - {\bm p}_3| $.
\end{description}
Finally, we set $h_3 = |{\bm p}_1 - {\bm p}_4|$. We implicitly assume that ${\bm p}_1$ and ${\bm p}_4$ belong to the same half-space. Additionally, note that ${h_1 \approx h_T}$.
\end{Cond}

\begin{note}
As an example, we {define} the matrices $A_{T}$ as 
\begin{align*}
\displaystyle
A_{T} := 
\begin{pmatrix}
\cos \theta  & - \sin \theta \\
 \sin \theta & \cos \theta
\end{pmatrix}, \quad 
{A}_{T} := 
\begin{pmatrix}
 \cos \theta  & - \sin \theta & 0\\
 \sin \theta & \cos \theta & 0 \\
 0 & 0 & 1 \\
\end{pmatrix},
\end{align*}
where $\theta$ denotes the angle. 
\end{note}

\begin{note}
None of the lengths of the edges of a simplex or the measures of the simplex {are} changed by the transformation, i.e.,
\begin{align}
\displaystyle
h_i \leq  h_{T}, \quad i=1,\ldots,d. \label{ineq331}
\end{align}	
\end{note}

\subsection{Additional notations and assumptions} \label{addinot}
\begin{figure}[htbp]
  \includegraphics[keepaspectratio, scale=0.45]{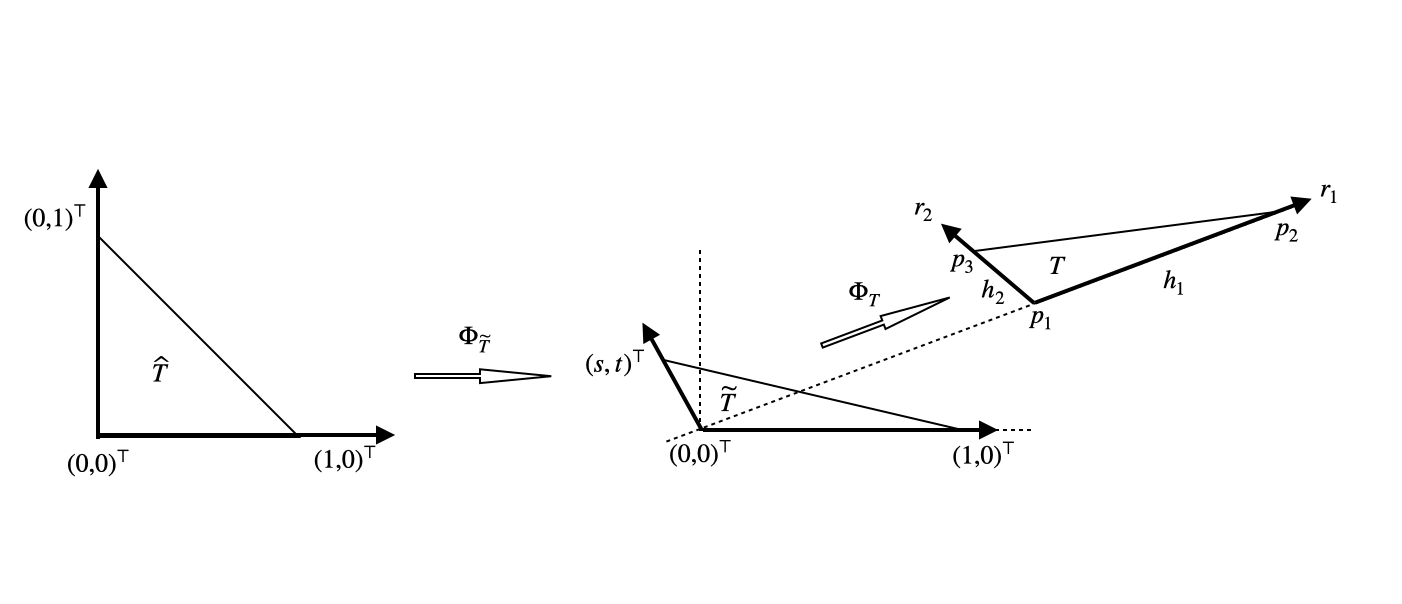}
\caption{Two-step affine mapping and vectors $r_i$, $i=1,2$}
\label{affine_2d}
\end{figure}

\begin{figure}[htbp]
  \begin{minipage}[b]{0.4\linewidth}
    \centering
    \includegraphics[keepaspectratio, scale=0.45]{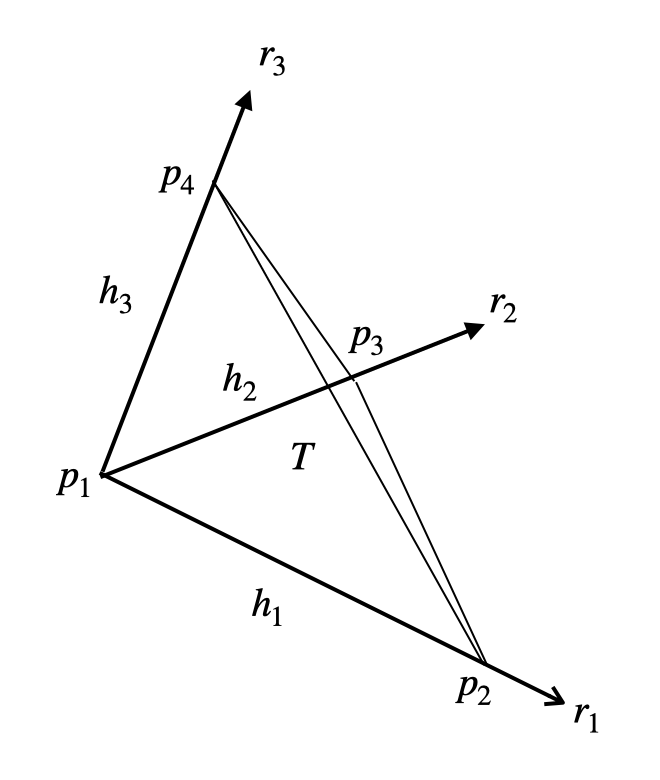}
    \caption{(Type \roman{sone}) Vectors $r_i$, $i=1,2,3$}
     \label{affine_3d_1}
  \end{minipage}
  \begin{minipage}[b]{0.4\linewidth}
    \centering
    \includegraphics[keepaspectratio, scale=0.45]{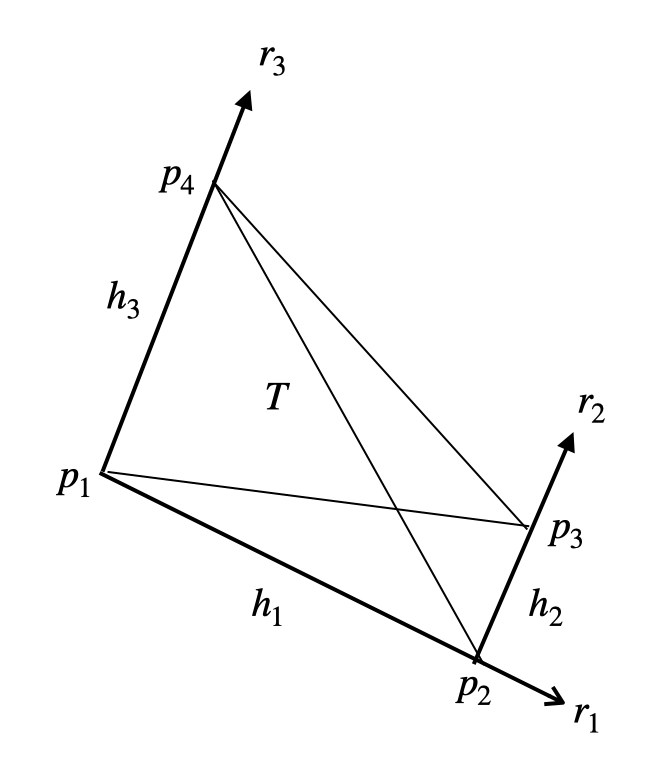}
    \caption{(Type \roman{stwo}) Vectors $r_i$, $i=1,2,3$}
     \label{affine_3d_2}
  \end{minipage}
\end{figure}

We define the vectors ${\bm{r}}_n \in \mathbb{R}^d$, $n=1,\ldots,d$ as follows: If $d=2$,
\begin{align*}
\displaystyle
{\bm r}_1 \coloneq \frac{{\bm p}_2 - {\bm p}_1}{|{\bm p}_2 - {\bm p}_1|}, \quad {\bm r}_2 \coloneq \frac{{\bm p}_3 - {\bm p}_1}{|{\bm p}_3 - {\bm p}_1|},
\end{align*}
see Fig. \ref{affine_2d}, and if $d=3$,
\begin{align*}
\displaystyle
&{\bm r}_1 \coloneq \frac{{\bm p}_2 - {\bm p}_1}{|{ \bm p}_2 - {\bm p}_1|}, \quad {\bm r}_3 \coloneq \frac{{\bm p}_4 - {\bm p}_1}{|{\bm p}_4 - {\bm p}_1|}, \quad
\begin{cases}
\displaystyle
{\bm r}_2 \coloneq \frac{{\bm p}_3 - {\bm p}_1}{|{\bm p}_3 - {\bm p}_1|}, \quad \text{for (Type \roman{sone})}, \\
\displaystyle
{\bm r}_2 \coloneq \frac{{\bm p}_3 - {\bm p}_2}{|{\bm p}_3 - {\bm p}_2|} \quad \text{for (Type \roman{stwo})},
\end{cases}
\end{align*}
see Fig \ref{affine_3d_1} for (Type \roman{sone}) and Fig \ref{affine_3d_2} for (Type \roman{stwo}). 

For a sufficiently smooth function $\varphi$ and a vector function ${\bm v} := (v_{1},\ldots,v_{d})^{\top}$, we define the directional derivative for $i \in \{ 1, \ldots,d \}$ as
\begin{align*}
\displaystyle
\frac{\partial \varphi}{\partial {{\bm r}_i}} &\coloneq ( {\bm r}_i \cdot  {\bm \nabla}_{x} ) \varphi = \sum_{i_0=1}^d ({\bm r}_i)_{i_0} \frac{\partial \varphi}{\partial x_{i_0}^{}}, \\
\frac{\partial {\bm v}}{\partial {\bm r}_i} &\coloneq \left(\frac{\partial v_{1}}{\partial {\bm r}_i}, \ldots, \frac{\partial v_{d}}{\partial {\bm r}_i} \right)^{\top} 
= ( ({\bm r}_i  \cdot {\bm \nabla}_{x}) v_{1}, \ldots, ({\bm r}_i  \cdot {\bm \nabla}_{x} ) v_{d} )^{\top}.
\end{align*}
For a multiindex $\beta = (\beta_1,\ldots,\beta_d) \in \mathbb{N}_0^d$, we use the notation
\begin{align}
\displaystyle
\partial^{\beta} \varphi \coloneq \frac{\partial^{|\beta|} \varphi}{\partial x_1^{\beta_1} \ldots \partial x_d^{\beta_d}}, \quad \partial^{\beta}_{r} \varphi \coloneq \frac{\partial^{|\beta|} \varphi}{\partial {\bm r}_1^{\beta_1} \ldots \partial {\bm r}_d^{\beta_d}}, \quad h^{\beta} \coloneq  h_{1}^{\beta_1} \cdots h_{d}^{\beta_d}. \label{partial=sym}
\end{align}
We note that $\partial^{\beta} \varphi \neq  \partial^{\beta}_{r} \varphi$.

 \begin{defi} \label{defi1}
 {We define the parameter $H_{{T}}$ as}
\begin{align*}
\displaystyle
H_{{T}} \coloneq \frac{\prod_{i=1}^d h_i}{|{T}|_d} h_{{T}}.
\end{align*}
\end{defi}
{This geometric condition was proposed in \cite{IshKobTsu21a,IshKobTsu21c} and is equivalent to the maximum angle condition \cite{IshKobSuzTsu21d}.}

\begin{assume} \label{neogeo=assume}
A family of meshes $\{ \mathbb{T}_h\}$ has a semi-regular property if there exists $\gamma_0 \> 0$ such that
\begin{align}
\displaystyle
\frac{H_{T}}{h_{T}} \leq \gamma_0 \quad \forall \mathbb{T}_h \in \{ \mathbb{T}_h \}, \quad \forall T \in \mathbb{T}_h. \label{NewGeo}
\end{align}
\end{assume}

\subsection{$L^2$-orthogonal projection}
For $T \in \mathbb{T}_h$, let $\Pi_{T}^0 : L^2(T) \to \mathbb{P}^0(T)$ be the $L^2$-orthogonal projection defined as
\begin{align*}
\displaystyle
\Pi_{T}^0 \varphi := \frac{1}{|T|_d} \int_{{T}} \varphi dx \quad \forall \varphi \in L^2(T).
\end{align*}
The following theorem provides an anisotropic error estimate for the projection $\Pi_{T}^0$. 
\begin{thr} \label{thr1}
For any ${\varphi} \in H^{1}({T})$,
\begin{align}
\displaystyle
\| \Pi_{T}^0 \varphi - \varphi \|_{L^2(T)} \leq c \sum_{i=1}^d h_i \left\| \frac{\partial \varphi}{\partial {\bm r}_i} \right\|_{L^{2}(T)}. \label{L2ortho}
\end{align}
\end{thr}

\begin{pf*}
This proof can be found in \cite[Theorem 2]{Ish24} and \cite[Theorem 2]{Ish24b}.
\qed
\end{pf*}

We also define the global interpolation $\Pi_h^0$ to space $P_{dc,h}^{0}$ as
\begin{align*}
\displaystyle
(\Pi_h^0 \varphi)|_{T} := \Pi_{{T}}^0 (\varphi|_{T}) \quad \forall T \in \mathbb{T}_h, \ \forall \varphi \in L^2(\Omega).
\end{align*}

\subsection{CR finite element interpolation operator}
For $T \in \mathbb{T}_h$, let the points $\{ P_{T,1}, \ldots, P_{T,d+1} \}$ be the {nodes} of the simplex $T \in \mathbb{T}_h$. Let $F_{T,i}$ be the face of $T$ opposite {to} $P_{T,i}$ for $i \in \{ 1, \ldots , d+1\}$. The CR interpolation operator $I_{T}^{CR} : H^{1}(T) \to \mathbb{P}^1(T)$ is defined as, for any $\varphi \in H^{1}(T)$, 
\begin{align*}
\displaystyle
I_{T}^{CR}: H^{1}(T) \ni \varphi \mapsto I_{T}^{CR} \varphi := \sum_{i=1}^{d+1} \left(  \frac{1}{| {F}_{T,i} |_{d-1}} \int_{{F}_{T,i}} {\varphi} d{s} \right) \theta_{T,i}^{CR} \in \mathbb{P}^1(T). 
\end{align*}
where $\theta_{T,i}^{CR}$ is the basis of the CR finite element (e.g. \cite[p. 82]{ErnGue21a}). We then present the estimates of the anisotropic CR interpolation error. 

\begin{thr} \label{DGCE=thr3}
For any ${\varphi} \in H^{2}({T})$,
\begin{align}
\displaystyle
|I_{T}^{CR} \varphi - \varphi |_{H^{1}({T})} &\leq c \sum_{i=1}^d {h}_i \left\| \frac{\partial }{\partial {\bm r}_i} \nabla \varphi \right \|_{L^2(T)^d}, \label{CRint}\\
\|I_{T}^{CR} \varphi - \varphi \|_{L^2(T)}
&\leq c  \sum_{|\varepsilon| = 2} h^{\varepsilon} \left\| \partial_{r}^{\varepsilon} \varphi  \right\|_{L^{2}(T)}. \label{CRint=L2}
\end{align}
Recall that $\partial_{r}^{\varepsilon}$ is defined in \eqref{partial=sym}.
\end{thr}

\begin{pf*}
The proof of \eqref{CRint} can be found in \cite[Theorem 3]{Ish24} and \cite[Theorem 3]{Ish24b}. The proof for \eqref{CRint=L2} can be found in \cite[Theorem 2]{Ish24c}.
\qed
\end{pf*}

We define the global interpolation operator ${I}_{h}^{CR}: H^{1}(\Omega) \to V_{h}^{CR}$ as
\begin{align*}
\displaystyle
({I}_{h}^{CR} \varphi )|_{T} = {I}_{T}^{CR} (\varphi |_{T}), \quad \forall T \in \mathbb{T}_h, \quad \forall \varphi \in H^{1}(\Omega),
\end{align*}	
where the space $V_{h}^{CR}$ is defined in \eqref{CRdef}.

\subsection{RT finite element interpolation operator} \label{RTsp}
For $T \in \mathbb{T}_h$, the local RT polynomial space is defined as
\begin{align}
\displaystyle
\mathbb{RT}^0(T) := \mathbb{P}^0(T)^d + {\bm x} \mathbb{P}^0(T), \quad {\bm x} \in \mathbb{R}^d. \label{RTsp}
\end{align}
Let $\mathcal{I}_{T}^{RT}: H^{1}(T)^d \to \mathbb{RT}^0(T)$ be the RT interpolation operator such that for any ${\bm v} \in H^{1}(T)^d$,
\begin{align}
\displaystyle
\mathcal{I}_{T}^{RT}: H^{1}(T)^d \ni {\bm v} \mapsto \mathcal{I}_{T}^{RT} {\bm v} := \sum_{i=1}^{d+1} \left(  \int_{{F}_{T,i}} {\bm v} \cdot {\bm n}_{T,i} d{s} \right) {\bm \theta}_{T,i}^{RT} \in \mathbb{RT}^0(T), \label{RTint}
\end{align}
where ${\bm \theta}_{T,i}^{RT}$ is the local shape {basis} function (e.g. \cite[p. 162]{ErnGue21a}) and ${\bm n}_{T,i}$ is a fixed unit normal to ${F}_{T,i}$. 

The following two theorems are divided into elements of (Type \roman{sone}) and (Type \roman{stwo}) in Section \ref{element=cond} when $d=3$.

\begin{thr} \label{DGRT=thr3}
Let $T$ be an element with Conditions \ref{cond1} or \ref{cond2} satisfying (Type \roman{sone}) in Section \ref{element=cond} when $d=3$. For any ${\bm v} = ({v}_1,\ldots,{v}_d)^T \in H^1(T)^d$,
\begin{align}
\displaystyle
\| \mathcal{I}_{T}^{RT} {\bm v} - {\bm v} \|_{L^2(T)^d} 
&\leq  c \left( \frac{H_{T}}{h_{T}} \sum_{i=1}^d h_i \left \|  \frac{\partial {\bm v}}{\partial {\bm r}_i} \right \|_{L^2(T)^d} +  h_{T} \| \div {\bm v} \|_{L^{2}({T})} \right). \label{RT5}
\end{align}
\end{thr}

\begin{pf*}
This proof is provided in \cite[Theorem 2]{Ish21}.
\qed
\end{pf*}

\begin{thr} \label{DGRT=thr4}
Let $d=3$. Let $T$ be an element with Condition \ref{cond2} that satisfies (Type \roman{stwo}) in Section \ref{element=cond}. For ${\bm v} = ({v}_1,v_2,{v}_3)^T \in H^1(T)^3$,
\begin{align}
\displaystyle
&\| \mathcal{I}_{T}^{RT} {\bm v} - {\bm v} \|_{L^2(T)^3} 
\leq c \frac{H_{T}}{h_{T}} \Biggl(  h_{T} |{\bm v}|_{H^1(T)^3} \Biggr). \label{RT6}
\end{align}
\end{thr}

\begin{pf*}
This proof is provided in \cite[Theorem 3]{Ish21}.
\qed
\end{pf*}

\begin{rem}
Below, we use the interpolation error estimate \eqref{RT5} in the ${\bm r}_i$-coordinate system for analysis.	
\end{rem}

The RT finite-element space is defined as follows:
\begin{align*}
\displaystyle
V^{RT}_{h} &:= \{ {\bm v_h} \in L^1(\Omega)^d: \  {\bm v_h}|_T \in \mathbb{RT}^0(T), \ \forall T \in \mathbb{T}_h, \  [\![ {\bm v_h} \cdot {\bm n} ]\!]_F = 0, \ \forall F \in \mathcal{F}_h^i \}.
\end{align*}
We define the following global RT interpolation $\mathcal{I}_{h}^{RT} : H^{1}(\Omega)^d \to V^{RT}_{h}$ as
\begin{align*}
\displaystyle
(\mathcal{I}_{h}^{RT} {\bm v} )|_{T} = \mathcal{I}_{T}^{RT} ({\bm v}|_{T}) \quad \forall T \in \mathbb{T}_h, \quad \forall {\bm v} \in H^{1}(\Omega)^d.
\end{align*}

\section{Discrete Poincar\'e inequality} \label{sec:5}

\setcounter{section}{4} \setcounter{equation}{0} 
This section presents the discrete Poincar\'e inequality on anisotropic meshes in Nitsche's method. Thus far, the following results have been obtained. A previous study \cite{Ish24} showed the discrete Poincar\'e inequality for the anisotropic weakly over-penalised symmetric interior penalty method. In \cite{Ish24c}, we presented the discrete Poincar\'e inequality for a hybrid anisotropic weakly over-penalised symmetric interior penalty method.  In \cite{Ish24d}, we introduced the discrete Sobolev inequality on anisotropic meshes for the CR finite element method. In \cite{Bre03,PieErn12}, the Poincar\'e inequality was proven under the shape-regularity condition. However, it is difficult to derive the inequality on anisotropic meshes. Herein, we present a proof using a dual problem. Therefore, we impose that $\Omega$ is convex to establish the aforementioned inequality.  

The following relation plays an important role in the Nitsche (or CR finite element) method for anisotropic meshes:

\begin{lem} \label{Nit=lem1}
For any ${\bm w} \in H^1(\Omega)^d$ and $\psi_h \in V_{h}^{CR} + H^1(\Omega)$,
\begin{align}
\displaystyle
&\int_{\Omega} \left( \mathcal{I}_h^{RT} {\bm w} \cdot {\bm \nabla}_h \psi_{h} + \div \mathcal{I}_h^{RT} {\bm w}  \psi_{h} \right) dx = \sum_{F \in \mathcal{F}_h^{\partial}} \int_{F} ( {\bm w} \cdot {\bm n}_F) \Pi_F^0 \psi_{h} ds.  \label{disPoi=1}
\end{align}
\end{lem}

\begin{pf*}
For any ${\bm w} \in H^1(\Omega)^d$ and $\psi_h \in V_{h}^{CR} + H^1(\Omega)$, using Green's formula and the fact $\mathcal{I}_h^{RT} {\bm w} \cdot {\bm n}_F \in \mathbb{P}^{0}(F)$ for any $F \in \mathcal{F}_h$ and $ [\![ \mathcal{I}_h^{RT} {\bm w} \cdot {\bm n}]\!]_F = 0$  for any $F \in \mathcal{F}_h^i$, we derive
\begin{align*}
\displaystyle
&\int_{\Omega} \left( \mathcal{I}_h^{RT} {\bm w} \cdot {\bm \nabla}_h \psi_{h} + \div \mathcal{I}_h^{RT} {\bm w}  \psi_{h} \right) dx
= \sum_{T \in \mathbb{T}_h} \int_{\partial T} (\mathcal{I}_h^{RT} {\bm w} \cdot {\bm n}_T) \psi_{h} ds  \notag \\
&\quad = \sum_{F \in \mathcal{F}_h^i} \int_{F} \left(  [\![ \mathcal{I}_h^{RT} {\bm w} \cdot {\bm n}]\!]_F \{\! \{ \psi_{h} \} \! \}_{\overline{\omega},F} + \{\! \{ \mathcal{I}_h^{RT} {\bm w}  \}\! \}_{\omega,F} \cdot {\bm n}_F [\![ \psi_{h} ]\!]_F \right) ds  \notag \\
&\quad + \sum_{F \in \mathcal{F}_h^{\partial}} \int_{F} (\mathcal{I}_h^{RT} {\bm w} \cdot {\bm n}_F) \psi_{h} ds \\
&\quad = \sum_{F \in \mathcal{F}_h^{\partial}} \int_{F} (\mathcal{I}_h^{RT} {\bm w} \cdot {\bm n}_F) \psi_{h} ds.
\end{align*}
Recall that the weighted and skew-weighted averages $ \{\! \{ \cdot \} \! \}_{{\omega},F}$ and  $ \{\! \{ \cdot \} \! \}_{\overline{\omega},F}$ were defined in Section \ref{sec=notation}. Using the properties of the projection $\Pi_F^0$ and $\mathcal{I}_h^{RT}$, we have
\begin{align*}
\displaystyle
\int_{\Omega} \left( \mathcal{I}_h^{RT} {\bm w} \cdot {\bm \nabla}_h \psi_{h} + \div \mathcal{I}_h^{RT} {\bm w}  \psi_{h} \right) dx
=  \sum_{F \in \mathcal{F}_h^{\partial}} \int_{F} ( {\bm w} \cdot {\bm n}_F) \Pi_F^0 \psi_{h} ds,
\end{align*}
which is the target inequality \eqref{disPoi=1}.
\qed
\end{pf*}

The right-hand term in \eqref{disPoi=1} is estimated as

\begin{lem} \label{Nit=lem2}
For any ${\bm w} \in H^1(\Omega)^d$ and $\psi_h \in V_{h}^{CR}$,
\begin{align}
\displaystyle
\left| \sum_{F \in \mathcal{F}_h^{\partial}} \int_{F} ( {\bm w} \cdot {\bm n}_F) \Pi_F^0 \psi_{h} ds \right| 
&\leq c |\psi_{h}|_1 \left(  h \| {\bm w} \|_{L^2(\Omega)^d} + h^{\frac{3}{2}} \| {\bm w} \|_{L^2(\Omega)^d}^{\frac{1}{2}} | {\bm w} |_{H^1(\Omega)^d}^{\frac{1}{2}} \right), \label{disPoi=2} \\
\left| \sum_{F \in \mathcal{F}_h^{\partial}} \int_{F} ( {\bm w} \cdot {\bm n}_F) \Pi_F^0 \psi_{h} ds \right|
&\leq  c |\psi_{h}|_{0}  \| {\bm w} \|_{H^1(\Omega)^d}. \label{disPoi=3}
\end{align}
\end{lem}

\begin{pf*}
Because $V_h^{CR} \subset  P_{dc,h}^{1}$, we can use the proofs in \cite[Lemmata 4 and 5]{Ish24}.
\qed
\end{pf*}

The following lemma provides the discrete Poincar\'e inequality:

\begin{lem}[Discrete Poincar\'e inequality] \label{Nit=lem3}
It is assumed that $\Omega$ is convex.  Let $\{ \mathbb{T}_h\}$ be a family of meshes with the semi-regular property (Assumption \ref{neogeo=assume}) and $h \leq 1$. Then, there exists a positive constant $C_{dc}^{P}$ independent of $h$ {but dependent on the maximum angle} such that
\begin{align}
\displaystyle
 \| \psi_{h} \|_{L^2(\Omega)} \leq C_{dc}^{P} | \psi_{h} |_{1,0} \quad \forall \psi_h \in V_{h}^{CR}. \label{disPoi=4}
\end{align}
\end{lem}

\begin{pf*}
For any $V_h^{CR} \subset  P_{dc,h}^{1}$, we can obtain the target inequality \eqref{disPoi=4} using \eqref{disPoi=1}, \eqref{disPoi=3} and the proof of \cite[Lemma 6]{Ish24}.
\qed
\end{pf*}

\begin{rem}
For the proof of Lemma \ref{Nit=lem3}, it is assumed that $\Omega$ is convex. Several methods have been proposed to avoid restrictive conditions (\cite{Bre03,PieErn12,GirLiRiv16,LasSul03}). However, these methods have been proven under the shape-regularity mesh condition. As described in \cite[Section 4.8]{Ish24d},  anisotropic meshes can be classified into two types: 
\begin{description}
  \item[(a)] those that include elements with large aspect ratios and those that do not satisfy the shape-regularity condition but satisfy the semi-regular condition (Assumption \ref{neogeo=assume});
  \item[(b)] those that include elements with large aspect ratios and whose partitions satisfy the shape-regularity condition.  
\end{description}
In (b), it may be possible to extend these new methods; however, it is difficult to apply these methods to the anisotropic meshes in (a). Further investigation of this issue remains a topic for future work. 
\end{rem}

\section{Stability and error estimates of the new scheme} \label{sec:6}
This section presents stability and error estimates of the scheme \eqref{ell_dis}.


As discussed in the introduction, the error \eqref{error=SymNit} in the symmetric Nitsche's method is derived by employing discrete stability \eqref{sta=SymNit} and Galerkin orthogonality \eqref{Galerkin=ortho}. The ultimate calculation of this error relies on the evaluation of the interpolation error.

Conversely, the scheme \eqref{ell_dis} lacks the inclusion of consistency terms, resulting in the absence of Galerkin orthogonality. Consequently, an alternative approach is required to derive error estimates. In this context, we propose a method utilising Lemma \ref{Nit=lem1}.

\setcounter{section}{5} \setcounter{equation}{0} 

\color{black}
\subsection{Stability}

\begin{thr}[Stability]
Assume that $\Omega$ is convex. Let $\{ \mathbb{T}_h\}$ be a family of meshes with the semiregular property (Assumption \ref{neogeo=assume}) and $h \leq 1$. Let $u_h  \in V_{h}^{CR}$ be a discrete solution to \eqref{ell_dis}. Then,
\begin{align}
\displaystyle
|u_h|_{1,1} \leq c \left(  \| f \|_{L^2(\Omega)} +  |g|_{1}  \right). \label{stability}
\end{align}
\end{thr}

\begin{pf*}
The H\"older, the Cauchy--Schwarz and the discrete Poincar\'e \eqref{disPoi=4} inequalities yield
\begin{align*}
\displaystyle
| \ell_h(u_h) |
&\leq  \| f \|_{L^2(\Omega)} \| u_h \|_{L^2(\Omega)} + \sum_{F \in \mathcal{F}_h^{\partial}} \kappa_{F(1)} \| \Pi_F^0  g \|_{L^2(F)} \| \Pi_F^0 u_h \|_{L^2(F)} \\
&\leq c  \| f \|_{L^2(\Omega)} |u_h|_{1,0} + |g|_{1} |u_h|_{1} \leq c \left(  \| f \|_{L^2(\Omega)} +  |g|_{1}  \right)  |u_h|_{1,1}
\end{align*}
if $\Omega$ is convex and $h \leq 1$. Here, we used the fact that $ | u_h |_{1,0} \leq  | u_h |_{1,1}$ for $h \leq 1$.
\qed
\end{pf*}

\subsection{Energy norm error estimate}
The starting point for the error analysis is the second Strang lemma \cite[P. 51]{ErnGue21b}.

\begin{lem} \label{second=starang}
It is assumed that  $\Omega$ is convex. Let $f \in L^2(\Omega)$ and $g \in H^{\frac{1}{2}}(\partial \Omega)$. Let $u \in H^1(\Omega)$ be the solution to \eqref{ell_eq}.  Let $u_h  \in V_{h}^{CR}$ be a discrete solution to \eqref{ell_dis}. Then,
\begin{align}
\displaystyle
| u - u_h |_{1,1} \leq c \inf_{v_ h \in V_h^{CR}} |u - v_h |_{1,1} + E_h(u), \label{po=err=2}
\end{align}
where
\begin{align}
\displaystyle
E_h(u) \coloneq \sup_{w_h \in V_h^{CR}} \frac{ | a_h(u , w_h) - \ell_h(w_h) |}{| w_h  |_{1,1}}. \label{po=err=3}
\end{align}
\end{lem}

\begin{pf*}
For any $v_h \in V_h^{CR}$, 
\begin{align*}
\displaystyle
| v_h - u_h |_{1,1}^2 
&= a_h(v_h - u_h , v_h - u_h) \\
&= a_h(v_h - u , v_h - u_h)  + a_h(u , v_h - u_h) - \ell_h(v_h - u_h) \\
&\leq c |v_h - u|_{1,1} |v_h - u_h|_{1,1} + | a_h(u , v_h - u_h) - \ell_h(v_h - u_h) |,
\end{align*}
which leads to
\begin{align*}
\displaystyle
| v_h - u_h |_{1,1}
&\leq c |u - v_h |_{1,1} + \frac{ | a_h(u , v_h - u_h) - \ell_h(v_h - u_h) |}{| v_h - u_h |_{1,1}} \\
&\leq  c |u - v_h |_{1,1} + E_h(u).
\end{align*}
Then,
\begin{align*}
\displaystyle
|u - u_h|_{1,1}
&\leq |u - v_h|_{1,1} + |v_h - u_h|_{1,1} \leq c  |u - v_h|_{1,1}  + E_h(u).
\end{align*}
Hence, the target inequality \eqref{po=err=2} holds true.
\qed
\end{pf*}

\begin{lem}[Best approximation] \label{best=approx}
We assume that $\Omega$ is convex. Let $f \in L^2(\Omega)$ and $g \in H^{\frac{1}{2}}(\partial \Omega)$. Let $u \in H^2(\Omega)$ be the solution to \eqref{ell_eq}. Then,
\begin{align}
\displaystyle
\inf_{v_ h \in V_h^{CR}} |u - v_h |_{1,1} \leq c \left( \sum_{i=1}^d \sum_{T \in \mathbb{T}_h} {h}_i^2 \left\| \frac{\partial }{\partial {\bm r}_i} {\bm \nabla} u \right \|_{L^2(T)^d}^2 \right)^{\frac{1}{2}}.  \label{po=err=4}
\end{align}
\end{lem}

\begin{pf*}
From the definitions of CR interpolation and $L^2$-projection, it holds that for any $F \in \mathcal{F}_h^{\partial}$,
\begin{align*}
\displaystyle
\Pi_F^0 (u - I_h^{CR} u) = \frac{1}{|F|_{d-1}} \int_F (u - I_h^{CR} u) ds = 0.
\end{align*}
Therefore, using \eqref{CRint}, we obtain
\begin{align*}
\displaystyle
\inf_{v_ h \in V_h^{CR}} |u - v_h |_{1,1}
&\leq |u - I_h^{CR} u |_{1,1} =  |u - I_h^{CR} u |_{H^1(\mathbb{T}_h)} \\
&\leq c \left( \sum_{i=1}^d \sum_{T \in \mathbb{T}_h} {h}_i^2 \left\| \frac{\partial }{\partial {\bm r}_i} {\bm \nabla} u \right \|_{L^2(T)^d}^2 \right)^{\frac{1}{2}},
\end{align*}
which is the target inequality \eqref{po=err=4}.
\qed
\end{pf*}

The essential part for the error estimate is the consistency error term \eqref{po=err=3}.

\begin{lem}[Asymptotic consistency] \label{asy=con}
It is assumed that $\Omega$ is convex. Let $f \in L^2(\Omega)$ and $g \in H^{\frac{1}{2}}(\partial \Omega)$. Let $u \in H^2(\Omega)$ be the solution to \eqref{ell_eq}. Let $\{ \mathbb{T}_h\}$ be a family of conformal meshes with semi-regular properties (Assumption \ref{neogeo=assume}). Let $T \in \mathbb{T}_h$ be the element with Conditions \ref{cond1} or \ref{cond2} satisfying (Type \roman{sone}) in Section \ref{element=cond} when $d=3$. Then,
\begin{align}
\displaystyle
E_h(u) &\leq c \left\{ \left( \sum_{i=1}^d \sum_{T \in \mathbb{T}_h} h_i^2 \left \| \frac{\partial}{\partial {\bm r}_i} {\bm \nabla} u \right \|_{L^2(T)^d}^2 \right)^{\frac{1}{2}} + h \| \varDelta u \|_{L^2(\Omega)} \right\} \notag \\
&\quad + c \left ( h |u|_{H^1(\Omega)} + h^{\frac{3}{2}} |u|_{H^1(\Omega)}^{\frac{1}{2}} \| \varDelta u\|_{L^2(\Omega)}^{\frac{1}{2}}  \right). \label{po=err=5}
\end{align}

\end{lem}

\begin{pf*}
Let $u \in H^2(\Omega)$. First, we have
\begin{align}
\displaystyle
 \div \mathcal{I}_{h}^{RT} {\bm \nabla} u &= \Pi_h^{0} \div {\bm \nabla} u =  \Pi_h^{0} \varDelta u.\label{po=err=6}
\end{align}
See \cite[Lemma 16.2]{ErnGue21a}. For any $w_h \in V_h^{CR}$, Setting ${\bm w} := {\bm \nabla} u$ in \eqref{disPoi=1} yields
\begin{align*}
\displaystyle
 &a_h(u , w_h) - \ell_h(w_h) \\
 &\quad = ( {\bm \nabla} u - \mathcal{I}_h^{RT} {\bm \nabla} u , {\bm \nabla}_h w_h  )_{\Omega} + ( \varDelta u - \Pi_h^0 \varDelta u , w_h)_{\Omega}
  + \sum_{F \in \mathcal{F}_h^{\partial}} \langle  {\bm \nabla} u \cdot {\bm n}_F , \Pi_F^0 w_{h} \rangle_F \\
 &\quad \eqcolon I_1 + I_2 + I_3.
\end{align*}
Using the H\"older inequality, the Cauchy--Schwarz inequality, and the RT interpolation error \eqref{RT5}, the term $I_1$ is estimated as
\begin{align*}
\displaystyle
|I_1|
&\leq c  \sum_{T \in \mathbb{T}_h} \| {\bm \nabla} u - \mathcal{I}_{h}^{RT} {\bm \nabla} u \|_{L^2(T)} | w_{h} |_{H^1(T)} \\
&\leq c \left\{ \left( \sum_{i=1}^d \sum_{T \in \mathbb{T}_h} h_i^2 \left \| \frac{\partial}{\partial {\bm r}_i} {\bm \nabla} u \right \|_{L^2(T)^d}^2 \right)^{\frac{1}{2}} + h \| \varDelta u \|_{L^2(\Omega)} \right\}  | w_{h} |_{1,1}.
\end{align*}
Using the H\"older inequality, the Cauchy--Schwarz inequality, the stability of $\Pi_h^0$, and the estimate \eqref{L2ortho}, the term $I_2$ is estimated as
\begin{align*}
\displaystyle
|I_2|
&= \left|  \int_{\Omega} \left( \varDelta u -  \Pi_h^0 \varDelta u \right) \left(  w_{h} - \Pi_h^0 w_{h}  \right) dx \right| \\
&\leq \sum_{T \in \mathbb{T}_h} \| \varDelta u -  \Pi_h^0 \varDelta u \|_{L^2(T)}  \| w_{h} - \Pi_h^0 w_{h} \|_{L^2(T)} \\
&\leq c h \| \varDelta u \|_{L^2(\Omega)}  | w_{h} |_{1,1}.
\end{align*}
Using inequality \eqref{disPoi=2}, the term $I_3$ is estimated as
\begin{align*}
\displaystyle
|I_3|
&\leq c \left(  h | u |_{H^1(\Omega)} + h^{ \frac{3}{2}} | u |_{H^1(\Omega)}^{\frac{1}{2}} \| \varDelta u \|_{L^2(\Omega)}^{\frac{1}{2}} \right)   | w_{h} |_{1,1}.
\end{align*}
Gathering the above inequalities, the target estimate \eqref{po=err=5} holds.
\qed
\end{pf*}

From {Lemma} \ref{second=starang}, \ref{best=approx}, and \ref{asy=con}, it is proven that the following energy norm error is estimated.

\begin{thr} \label{energy=thr}
It is assumed that $\Omega$ is convex. Let $f \in L^2(\Omega)$ and $g \in H^{\frac{1}{2}}(\partial \Omega)$. Let $u \in H^2(\Omega)$ be the solution to \eqref{ell_eq}. Let $\{ \mathbb{T}_h\}$ be a family of conformal meshes with semi-regular properties (Assumption \ref{neogeo=assume}). Let $T \in \mathbb{T}_h$ be the element with Conditions \ref{cond1} or \ref{cond2} satisfying (Type \roman{sone}) in Section \ref{element=cond} when $d=3$. Let $u_h  \in V_{h}^{CR}$ be a discrete solution of \eqref{ell_dis}. Then,
\begin{align}
\displaystyle
|u-u_h|_{1,1} &\leq c \left\{ \left( \sum_{i=1}^d \sum_{T \in \mathbb{T}_h} h_i^2 \left \| \frac{\partial}{\partial {\bm r}_i} {\bm \nabla} u \right \|_{L^2(T)^d}^2 \right)^{\frac{1}{2}} + h \| \varDelta u \|_{L^2(\Omega)} \right\} \notag \\
&\quad + c \left ( h |u|_{H^1(\Omega)} + h^{\frac{3}{2}} |u|_{H^1(\Omega)}^{\frac{1}{2}} \| \varDelta u\|_{L^2(\Omega)}^{\frac{1}{2}}  \right). \label{po=err=7}
\end{align}
\end{thr}

\subsection{$L^2$ norm error estimate}
This section presents the $L^2$-error estimate for the Nitsche method of the Poisson equation.

\begin{thr} \label{L2=thr}
We assume that $\Omega$ is convex. Let $f \in L^2(\Omega)$ and $g \in H^{\frac{1}{2}}(\partial \Omega)$. Let $u \in H^2(\Omega)$ be the solution to \eqref{ell_eq}. Let $\{ \mathbb{T}_h\}$ be a family of conformal meshes with semi-regular properties (Assumption \ref{neogeo=assume}). Let $T \in \mathbb{T}_h$ be the element with Conditions \ref{cond1} or \ref{cond2} satisfying (Type \roman{sone}) in Section \ref{element=cond} when $d=3$. Let $u_h  \in V_{h}^{CR}$ be a discrete solution of \eqref{ell_dis}. Then,
\begin{align}
\displaystyle
\| u - u_h \|_{L^2(\Omega)}
&\leq c h \left\{ \left( \sum_{i=1}^d \sum_{T \in \mathbb{T}_h} h_i^2 \left \| \frac{\partial}{\partial {\bm r}_i} {\bm \nabla} u \right \|_{L^2(T)^d}^2 \right)^{\frac{1}{2}} + h \| \varDelta u \|_{L^2(\Omega)} \right\} \notag \\
&\quad + c h  \left ( h |u|_{H^1(\Omega)} + h^{\frac{3}{2}} |u|_{H^1(\Omega)}^{\frac{1}{2}} \| \varDelta u\|_{L^2(\Omega)}^{\frac{1}{2}}  \right) \notag \\
&\quad + c \left( \sum_{|\varepsilon| = 2}  \sum_{T \in \mathbb{T}_h} h^{2 \varepsilon} \left\| \partial_{r}^{ \varepsilon} u  \right\|_{L^{2}(T)}^2 \right)^{\frac{1}{2}}. \label{po=err=8}
\end{align}
\end{thr}

\begin{pf*}
We set $e:=  u - u_h$.  Let $z \in H_0^1(\Omega) \cap H^2(\Omega)$ {satisfying}
\begin{align}
\displaystyle
- \varDelta z = e \quad \text{in $\Omega$}, \quad z = 0 \quad \text{on $\partial \Omega$}. \label{po=err=9}
\end{align}
Note that $|z|_{H^1(\Omega)} \leq c \| e \|_{L^2(\Omega)}$ and $|z|_{H^2(\Omega)} \leq \| e \|_{L^2(\Omega)}$.  Let $z_h \in V_{h}^{CR}$ satisfy
\begin{align}
\displaystyle
a_h(\varphi_h , z_h) = (\varphi_h , e)_{\Omega} \quad \forall \varphi_h \in V_{h}^{CR}. \label{po=err=10}
\end{align}
Then,
\begin{align*}
\displaystyle
\| e \|^2_{L^2(\Omega)}
&= (u - u_h , e)_{\Omega} = a_h(u,z) - \int_{\partial  \Omega} {\bm \nabla} z \cdot {\bm n} u ds - a_h(u_h,z_h) \\
&= a_h(u-u_h,z-z_h) + a_h(u_h,z-z_h)+a_h(u-u_h,z_h) - \int_{\partial  \Omega} {\bm \nabla} z \cdot {\bm n} u ds \\
&=  a_h(u-u_h,z-z_h) \\
&\quad + a_h(u-u_h,z_h - I_h^{CR} z) + a_h(u-u_h,I_h^{CR} z) \\
&\quad + a_h(u_h - I_h^{CR} u, z - z_h) + a_h(I_h^{CR} u, z - z_h) - \int_{\partial  \Omega} {\bm \nabla} z \cdot {\bm n} u ds \\
&\eqcolon J_1 + J_2 + J_3 + J_4 + J_5 + J_6.
\end{align*}
From Theorems \ref{DGCE=thr3} and \ref{energy=thr} with $g \equiv 0$, we obtain
\begin{subequations} \label{po=err=11}
\begin{align}
\displaystyle
 | z- z_h |_{1,1}
 &\leq c h \| e \|_{L^2(\Omega)}, \label{po=err=11a} \\
 | z- I_{h}^{CR} z |_{H^1(\mathbb{T}_h)}
  &\leq c  h \| e \|_{L^2(\Omega)}, \label{po=err=11b} \\
 \| z- I_{h}^{CR} z \|_{L^2(\Omega)}
  &\leq c  h^2 \| e \|_{L^2(\Omega)}. \label{po=err=11c}
\end{align}
\end{subequations}
Using \eqref{po=err=7} and \eqref{po=err=11a}, $J_1$ can be estimated as
\begin{align*}
\displaystyle
|J_1| 
&\leq c |u - u_h|_{1,1} |z - z_h|_{1,1} \\
&\leq c  h \| e \|_{L^2(\Omega)} \left\{ \left( \sum_{i=1}^d \sum_{T \in \mathbb{T}_h} h_i^2 \left \| \frac{\partial}{\partial {\bm r}_i} {\bm \nabla} u \right \|_{L^2(T)^d}^2 \right)^{\frac{1}{2}} + h \| \varDelta u \|_{L^2(\Omega)} \right\} \notag \\
&\quad + c  h \| e \|_{L^2(\Omega)} \left ( h |u|_{H^1(\Omega)} + h^{\frac{3}{2}} |u|_{H^1(\Omega)}^{\frac{1}{2}} \| \varDelta u\|_{L^2(\Omega)}^{\frac{1}{2}}  \right).
\end{align*}
Using \eqref{po=err=7}, \eqref{po=err=11a} and \eqref{po=err=11b}, $J_2$ can be estimated as
\begin{align*}
\displaystyle
|J_2|
&\leq c | u- u_h |_{1,1} \left(  | z_h - z |_{1,1} + | z - I_h^{CR} z|_{H^1(\mathbb{T}_h)} \right) \\
&\leq c h \| e \|_{L^2(\Omega)}  \left\{ \left( \sum_{i=1}^d \sum_{T \in \mathbb{T}_h} h_i^2 \left \| \frac{\partial}{\partial {\bm r}_i} {\bm \nabla} u \right \|_{L^2(T)^d}^2 \right)^{\frac{1}{2}} + h \| \varDelta u \|_{L^2(\Omega)} \right\} \notag \\
&\quad + c h \| e \|_{L^2(\Omega)} \left ( h |u|_{H^1(\Omega)} + h^{\frac{3}{2}} |u|_{H^1(\Omega)}^{\frac{1}{2}} \| \varDelta u\|_{L^2(\Omega)}^{\frac{1}{2}}  \right).
\end{align*}
By using an analogous argument,
\begin{align*}
\displaystyle
|J_4|
&\leq c  \left(  | u_h - u |_{1,1} + | u - I_h^{CR} u|_{H^1(\mathbb{T}_h)} \right)  | z- z_h |_{1,1} \\
&\leq c h \| e \|_{L^2(\Omega)}  \left\{ \left( \sum_{i=1}^d \sum_{T \in \mathbb{T}_h} h_i^2 \left \| \frac{\partial}{\partial {\bm r}_i} {\bm \nabla} u \right \|_{L^2(T)^d}^2 \right)^{\frac{1}{2}} + h \| \varDelta u \|_{L^2(\Omega)} \right\} \notag \\
&\quad + c h \| e \|_{L^2(\Omega)} \left ( h |u|_{H^1(\Omega)} + h^{\frac{3}{2}} |u|_{H^1(\Omega)}^{\frac{1}{2}} \| \varDelta u\|_{L^2(\Omega)}^{\frac{1}{2}}  \right).
\end{align*}
Using Lemma \ref{Nit=lem1} and the definition of CR interpolation,
\begin{align}
\displaystyle
&\int_{\Omega} \left( \mathcal{I}_h^{RT} {\bm \nabla} z \cdot {\bm \nabla}_h (I_{h}^{CR} u - u ) + \div \mathcal{I}_h^{RT} ({\bm \nabla} z)  (I_{h}^{CR} u - u ) \right) dx \notag\\
&\quad =  \sum_{F \in \mathcal{F}_h^{\partial}} \int_{F} ( {\bm \nabla} z \cdot {\bm n}_F) \Pi_F^0(I_{h}^{CR} u - u) ds = 0. \label{po=err=12}
\end{align}
From \eqref{po=err=9},
\begin{align*}
\displaystyle
 - (u,e)_{\Omega}
 &=  (u,\varDelta z)_{\Omega} = \int_{\partial \Omega} {\bm \nabla} z \cdot {\bm n} u ds -  ({\bm \nabla} u, {\bm \nabla} z)_{\Omega}.
\end{align*}
Thus, using \eqref{po=err=10} and \eqref{po=err=12}, we obtain
\begin{align*}
\displaystyle
J_5+J_6
&=  a_h(I_h^{CR} u, z ) - a_h(I_h^{CR} u, z_h) - \int_{\partial  \Omega} {\bm \nabla} z \cdot {\bm n} u ds \\
&=  a_h(I_h^{CR} u - u, z ) + a_h(u,z) - (I_h^{CR} u - u , e)_{\Omega} - (u,e)_{\Omega} - \int_{\partial  \Omega} {\bm \nabla} z \cdot {\bm n} u ds \\
&= a_h(I_h^{CR} u - u, z ) - (I_h^{CR} u - u , e)_{\Omega} \\
&=  \int_{\Omega} ({\bm \nabla} z -  \mathcal{I}_h^{RT} {\bm \nabla} z  ) \cdot {\bm \nabla}_h ( I_{h}^{CR} u - u ) dx + \int_{\Omega} ( \varDelta z - \Pi_h^0 \varDelta z ) (I_{h}^{CR} u - u )dx.
\end{align*}
Using \eqref{CRint}, \eqref{CRint=L2}, \eqref{RT5}, and the stability of the $L^2$-projection, we have
\begin{align*}
\displaystyle
|J_5 + J_6| 
&= \left| \int_{\Omega} ({\bm \nabla} z -  \mathcal{I}_h^{RT} {\bm \nabla} z  ) \cdot {\bm \nabla}_h ( I_{h}^{CR} u - u ) dx + \int_{\Omega} ( \varDelta z - \Pi_h^0 \varDelta z ) (I_{h}^{CR} u - u )dx \right| \\
&\leq  c h  \| e \|_{L^2(\Omega)}  \left(  \sum_{i=1}^d  \sum_{T \in \mathbb{T}_h} h_i^2 \left \| \frac{\partial}{\partial {\bm r}_i} {\bm \nabla} u \right \|_{L^2(T)^d}^2 \right)^{\frac{1}{2}} \\
&\quad+ c \| e \|_{L^2(\Omega)} \left( \sum_{|\varepsilon| = 2}  \sum_{T \in \mathbb{T}_h} h^{2 \varepsilon} \left\| \partial_{r}^{ \varepsilon} u  \right\|_{L^{2}(T)}^2 \right)^{\frac{1}{2}}.
\end{align*}
Using an analogous argument:
\begin{align*}
\displaystyle
J_3
&=  a_h(u,I_h^{CR} z) -  a_h(u_h,I_h^{CR} z) \\
&=  \int_{\Omega} ({\bm \nabla} u -  \mathcal{I}_h^{RT} {\bm \nabla} u  ) \cdot {\bm \nabla}_h ( I_{h}^{CR} z - z ) dx + \int_{\Omega} ( \varDelta u - \Pi_h^0 \varDelta u ) (I_{h}^{CR} z - z )dx,
\end{align*}
which leads to
\begin{align*}
\displaystyle
|J_3|
&\leq  c h  \| e \|_{L^2(\Omega)} \left\{ \left( \sum_{i=1}^d \sum_{T \in \mathbb{T}_h} h_i^2 \left \| \frac{\partial}{\partial {\bm r}_i} {\bm \nabla} u \right \|_{L^2(T)^d}^2 \right)^{\frac{1}{2}} + h \| \varDelta u \|_{L^2(\Omega)} \right\} \\
&\quad + c h^2  \| e \|_{L^2(\Omega)} \| \varDelta u \|_{L^2(\Omega)}.
\end{align*}
Combining the above inequalities yields the target inequality \eqref{po=err=8}.
\qed
\end{pf*}

\section{Pressure robust scheme for the Stokes problem} \label{sec:7}

\setcounter{section}{6} \setcounter{equation}{0} 
This section analyses the Nitsche method for the Stokes equations on anisotropic meshes.

\subsection{Continuous problem}
The (scaled) Stokes problem is to find $({\bm u},p): \Omega \to \mathbb{R}^d \times \mathbb{R}$ such that
\begin{align}
\displaystyle
- \nu \varDelta {\bm u} + {\bm \nabla} p = {\bm f} \quad \text{in $\Omega$}, \quad \div {\bm u}  = 0 \quad \text{in $\Omega$}, \quad {\bm u} = {\bm g} \quad \text{on $\partial \Omega$}, \label{stokes1}
\end{align}
where $\nu$ is a nonnegative parameter with $\nu \leq 1$, ${\bm f} \in L^2(\Omega)^d$ is a given function and ${\bm g} \in H^{\frac{1}{2}}(\partial \Omega)^d$ satisfying $\int_{\partial \Omega} {\bm g} \cdot {\bm n} ds = 0$. We consider the following function spaces:
\begin{align*}
\displaystyle
W({\bm g}) &\coloneq \left \{  {\bm v} \in H^1(\Omega)^d: \ {\bm v} = {\bm g} \text{ on } \partial \Omega \right \}, \\
Q &\coloneq L^2_0(\Omega) \coloneq \left \{ q \in L^2(\Omega); \ \int_{\Omega} q dx = 0 \right \},
\end{align*}
with norms:
\begin{align*}
\displaystyle
| \cdot |_{W(0)} := | \cdot |_{H^1(\Omega)^d}, \quad \| \cdot \|_Q := \| \cdot \|_{L^2(\Omega)}.
\end{align*}
The variational formulation for the Stokes equation \eqref{stokes1} is as follows: For any ${\bm f} \in L^2(\Omega)^d$ and ${\bm g} \in H^{\frac{1}{2}}(\partial \Omega)^d$, find $({\bm u},p) \in W({\bm g}) \times Q$ such that
\begin{subequations} \label{stokes2}
\begin{align}
\displaystyle
\nu A({\bm u},{\bm v}) + B({\bm v} , p) &= ({\bm f} , {\bm v})_{\Omega} \quad \forall {\bm v} \in W({\bm 0}), \label{stokes2a} \\
B({\bm u} , q) &= 0 \quad \forall q \in Q, \label{stokes2b}
\end{align}
\end{subequations}
where $A: H^1(\Omega)^d \times H^1(\Omega)^d \to \mathbb{R}$ and $B:H^1(\Omega)^d \times L^2(\Omega) \to \mathbb{R}$ respectively denote bilinear forms defined by
\begin{align*}
\displaystyle
A({\bm u},{\bm v}) \coloneq \int_{\Omega} {\bm \nabla} {\bm u} : {\bm \nabla} {\bm v} dx = \sum_{i=1}^d (  {\bm \nabla} u_i , {\bm \nabla} v_i  )_{\Omega}, \quad   B({\bm v}, q) \coloneq - ( \div {\bm v} ,  q )_{\Omega}. 
\end{align*}
Here, the colon denotes the scalar product of tensors. The continuous inf-sup inequality
\begin{align}
\displaystyle
\inf_{q \in Q} \sup_{{\bm v} \in W({\bm 0})} \frac{B({\bm v} , q)}{ | {\bm v} |_{W({\bm 0})} \| q \|_{Q} }\geq \beta \>0 \label{infsup=con}
\end{align}
has been shown to hold. Proofs can be found in Ern and Guermond \cite[Lemma 53.9]{ErnGue21b}, Girault and Raviart \cite[Lemma 4.1]{GirRav86}, and John \cite[Theorem 3.46]{Joh16}.

We set 
\begin{align*}
\displaystyle
H_*^1(\Omega) := H^1(\Omega) \cap L^2_0(\Omega), \quad \mathcal{H}
 := \{ {\bm v} \in L^2(\Omega)^d: \ \div {\bm v} = 0, \ {\bm v}|_{\partial \Omega} \cdot {\bm n} = 0 \},
 \end{align*}
 where $\div {\bm v} = 0$ and ${\bm v}|_{\partial \Omega} \cdot {\bm n} = 0$ indicate that $\int_{\Omega} ({\bm v} \cdot {\bm \nabla} ) q dx = 0$ for any $q \in H_*^1(\Omega) $. Then, the following $L^2$-orthogonal decomposition holds:
\begin{align}
\displaystyle
L^2(\Omega)^d = \mathcal{H} \oplus \nabla (H_*^1(\Omega) ). \label{helmdeco}
\end{align}
See \cite[Lemma 74.1]{ErnGue21c}. The $L^2$-orthogonal projection $P_{\mathcal{H}}: L^2(\Omega)^d \to \mathcal{H}$ resulting from this decomposition is often called \textit{Leray projection}.

From \cite[Lemma 2.2]{GirRav86}, there exists a function ${\bm u}_0 \in H^1(\Omega)^d$ such that
\begin{align*}
\displaystyle
\div {\bm u}_0 = 0 \quad \text{in $\Omega$}, \quad {\bm u}_0 = {\bm g} \quad \text{on $\partial \Omega$}.
\end{align*}
The function ${\bm w \coloneq \bm u- \bm{u}_0}$ satisfies the Stokes problem with the homogeneous Dirichlet boundary condition and a modified right-hand side:
\begin{align*}
\displaystyle
\ell({\bm v}) \coloneq ({\bm f , \bm v})_{\Omega}  - \nu A({\bm u_0,\bm v}) \quad \forall {\bm v} \in W({\bm 0}).
\end{align*}
We consider the following problem: Find $({\bm w},p) \in W({\bm 0}) \times Q$ such that
\begin{subequations} \label{stokes7}
\begin{align}
\displaystyle
\nu A({\bm w,\bm v}) + B({\bm v} , p) &= \ell({\bm v}) \quad \forall {\bm v} \in W({\bm 0}), \label{stokes7a} \\
B({\bm w} , q) &= 0 \quad \forall q \in Q, \label{stokes7b}
\end{align}
\end{subequations}

\begin{thr}
For any ${\bm f} \in L^2(\Omega)^d$, there exists a unique pair $({\bm w},p) \in W({\bm 0}) \times Q$ that satisfies \eqref{stokes1}.
\begin{align}
\displaystyle
| {\bm w} |_{W({\bm 0})} &\leq \frac{C_P}{\nu}\| P_{\mathcal{H}} ({\bm f}) \|_{L^2(\Omega)^d} +  | {\bm u}_0 |_{W({\bm 0})}, \label{stokes8} \\
\| p \|_{L^2(\Omega)} &\leq \frac{1+C_P}{\beta} \| {\bm f} \|_{L^2(\Omega)^d} + \frac{2 \nu}{\beta}  | {\bm u}_0 |_{W({\bm 0})}.  \label{stokes9}
\end{align}
\end{thr}

\begin{pf*}
The proof of the uniqueness of the solution of the Stokes equation \eqref{stokes7} can be found in \cite[Theorem 5.1]{GirRav86}.  Because ${\bm w}$ is divergence-free and vanishes at the boundary, we have ${\bm w} \in \mathcal{H}$, and
\begin{align*}
\displaystyle
 \int_{\Omega} {\bm f \cdot \bm w} dx = \int_{\Omega} P_{\mathcal{H}} ({\bm f}) \cdot {\bm w} dx.
\end{align*}
Setting ${\bm v := \bm w}$ in \eqref{stokes7a} and $q := p$ in \eqref{stokes7b} yields
\begin{align*}
\displaystyle
 \nu | {\bm w} |^2_{W({\bm 0})} &= \int_{\Omega} P_{\mathcal{H}} ({\bm f}) \cdot {\bm w} dx - \nu A({\bm u_0,\bm w}) \\
&\leq \| P_{\mathcal{H}} ({\bm f}) \|_{L^2(\Omega)^d} \| {\bm w} \|_{L^2(\Omega)^d} + \nu  | {\bm u}_0 |_{W({\bm 0})} | {\bm w} |_{W({\bm 0})}\\
&\leq C_P \| P_{\mathcal{H}} ({\bm f}) \|_{L^2(\Omega)^d} | {\bm w} |_{W({\bm 0})} + \nu  | {\bm u}_0 |_{W({\bm 0})} | {\bm w} |_{W({\bm 0})},
\end{align*}
which leads to \eqref{stokes8}. For the estimate of  the pressure, using the inf-sup condition \eqref{infsup=con}, the equation \eqref{stokes7a}, the H\"older's inequality, and \eqref{stokes8} yields
\begin{align*}
\displaystyle
\beta \| p \|_Q 
&\leq \sup_{\bm v \in W(\bm 0)} \frac{|B(\bm v,p)|}{|\bm v|_{W(\bm 0)}}
= \sup_{\bm v \in W(\bm 0)} \frac{|(\bm f,\bm v)_{\Omega} - \nu A(\bm u_0,\bm v) - \nu A(\bm w.\bm v)|}{|\bm v|_{W(\bm 0)}} \\
&\leq  \sup_{\bm v \in W(\bm 0)} \frac{ \| \bm f \|_{L^2(\Omega)^d} |\bm v|_{W(\bm 0)} + \nu |\bm u_0|_{W(\bm 0)} |\bm v|_{W(\bm 0)} + \nu |\bm w|_{W(\bm 0)} |\bm v|_{W(\bm 0)} }{|\bm v|_{W(\bm 0)}} \\
&=  \| \bm f \|_{L^2(\Omega)^d} + \nu |\bm u_0|_{W(\bm 0)} + \nu |\bm w|_{W(\bm 0)} \leq (1+C_P) \| \bm f \|_{L^2(\Omega)^d} + 2 \nu  | \bm u_0 |_{W(\bm 0)},
\end{align*}
which leads to \eqref{stokes9}. Here, we used that $ \| P_{\mathcal{H}} (\bm f) \|_{L^2(\Omega)^d} \leq  \| \bm f \|_{L^2(\Omega)^d} $.
\qed
\end{pf*}

\color{black}
\subsection{New scheme} \label{new=stokes}
Let $(W_{h}^{CR},Q_h^{0})$ be a pair of finite element spaces defined by
\begin{align}
\displaystyle
W_{h}^{CR} := (V_{h}^{CR})^d, \quad Q_h^{0} = P_{dc,h}^{0} \cap Q.\label{dis=stokes1}
\end{align}
The RT finite element space is defined as follows:
\begin{align*}
\displaystyle
V^{RT}_{h0} &:= \{ {\bm v}_h \in V^{RT}_{h}: \ ({\bm v}_h \cdot {\bm n})|_F = 0, \ \forall F \in \mathcal{F}_h^{\partial} \}.
\end{align*}
We define the following global RT interpolation $\mathcal{I}_{h0}^{RT} : W_h^{CR} \to V^{RT}_{h0}$ as
\begin{align*}
\displaystyle
(\mathcal{I}_{h0}^{RT} {\bm v}_h )|_{T} = \mathcal{I}_{T}^{RT} ({\bm v}_h|_{T}) \quad \forall T \in \mathbb{T}_h, \quad \forall {\bm v}_h \in  W_h^{CR}.
\end{align*}

For all ${\bm f} \coloneq (f_1,\ldots,f_d)^{\top} \in L^2(\Omega)^d$ and ${\bm g} \coloneq (g_1,\ldots,g_d)^{\top} \in H^{\frac{1}{2}}(\partial \Omega)^d$, we consider the Nitsche method for the Stokes equation \eqref{stokes1} as follows: We aim to find $({\bm u}_h,p_h) \in W_{h}^{CR} \times Q_h^{0}$ such that
\begin{subequations} \label{dis=stokes2}
\begin{align}
\displaystyle
\nu A_h(\bm u_h,\bm v_h) + B_h(\bm v_h , p_h) &= L_h(\bm v_h) \quad \forall \bm v_h \in W_{h}^{CR}, \label{dis=stokes2a} \\
B_h(\bm u_h , q_h) &= 0 \quad \forall q_h \in Q_h^{0}, \label{dis=stokes2b}
\end{align}
\end{subequations}
where $A_h: (H^1(\Omega)^d + W_{h}^{CR}) \times (H^1(\Omega)^d + W_{h}^{CR}) \to \mathbb{R}$, $B_h: (H^1(\Omega)^d+W_{h}^{CR}) \times Q_h^{0} \to \mathbb{R}$ and $L_h: H^1(\Omega)^d + W_h^{CR} \to \mathbb{R}$ are respectively defined as
\begin{align*}
\displaystyle
A_h(\bm u_h,\bm v_h) &\coloneq  \sum_{i=1}^d  \left\{ a_h^i(\bm u_{h},\bm v_{h}) +  s_h^i(\bm u_{h},\bm v_{h}) \right\}, \\
a_h^i(\bm u_{h},\bm v_{h}) &\coloneq  (\bm \nabla_h u_{h,i} , \bm \nabla_h v_{h,i})_{\Omega}, \\
s_h^i(\bm u_{h},\bm v_{h}) &\coloneq \eta \sum_{F \in \mathcal{F}_h^{\partial}} \kappa_{F(1)} \langle \Pi_F^0  u_{h,i}   , \Pi_F^0 v_{h,i} \rangle_F,\\
B_h(\bm v_h , q_h)& \coloneq - (\divh \bm v_h, q_h)_{\Omega}, \\
L_h(\bm v_h) &\coloneq \sum_{i=1}^d \ell^i_h(\bm v_h), \\
\ell^i_h(\bm v_h) &\coloneq ( f_i , (\mathcal{I}_{h0}^{RT} \bm v_h)_{i} )_{\Omega} + \nu \eta  \sum_{F \in \mathcal{F}_h^{\partial}} \kappa_{F(1)} \langle \Pi_F^0  g_{i}, \Pi_F^0 v_{h,i} \rangle_F,
\end{align*}
\color{black}
where $\eta$ is a positive parameter, $\{ u_{h,i}\}_{i=1}^d$, $\{ v_{h,i}\}_{i=1}^d$, $\{ f_{i}\}_{i=1}^d$ and $\{ (\mathcal{I}_{h0}^{RT} {\bm v}_h)_{i} \}_{i=1}^d$ denote the Cartesian components of ${\bm u}_h$, ${\bm v}_h$, ${\bm f}$ and $\mathcal{I}_{h0}^{RT} {\bm v}_h$. Here, for all $T \in \mathbb{T}_h$, we define the broken $H(\div;T)$ space as
\begin{align*}
\displaystyle
H(\div;\mathbb{T}_h) \coloneq \left \{ {\bm v} \in L^2(\Omega)^d; \ {\bm v} |_T \in H(\div;T) \ \forall T \in \mathbb{T}_h  \right\},
\end{align*}
and the broken divergence operator $\divh : H(\div;\mathbb{T}_h) \to L^2(\Omega)$ such that for all ${\bm v} \in H(\div;\mathbb{T}_h)$,
\begin{align*}
\displaystyle
(\divh {\bm v})|_T \coloneq \div ({\bm v} |_T) \quad \forall T \in \mathbb{T}_h.
\end{align*}

We define norms as follows:
\begin{align*}
\displaystyle
|{\bm w}_h|_{W_{h}^{CR}} := \left ( \sum_{i=1}^d |w_{h,i}|^2_{1,1} \right)^{\frac{1}{2}}, \quad \| q_h \|_{Q_h^{0}} := \| q_h \|_{L^2(\Omega)}
\end{align*}
for any ${\bm w}_h = (w_{h,1},\ldots,w_{h,d})^{\top}  \in W_{h}^{CR}$ and $q_h \in Q_h^{0}$.

We define a discrete divergence-free weak subspace as follows:
\begin{align*}
\displaystyle
W_{h,\div}^{CR} := \{ {\bm w}_h \in W_h^{CR}: \ B_h({\bm w}_h,q_h) = 0 \ \forall q_h \in Q_h \}.
\end{align*}
Then, the reduced problem of \eqref{dis=stokes2} is as follows. Find ${\bm u}_h \in W_{h,\div}^{CR}$ such that
\begin{align}
\displaystyle
\nu A_h({\bm u_h,\bm v_h})  &= L_h({\bm v}_h) \quad \forall {\bm v}_h \in W_{h,\div}^{CR}. \label{dis=stokes3}
\end{align}

\begin{rem} \label{dis=coe=rem}
The following discrete coercivity holds for any $\nu \> 0$ and $\eta \>0$: Indeed, for any ${\bm v}_h \in W_{h}^{CR}$,
\begin{align}
\displaystyle
\nu A_h({\bm v_h,\bm v_h})
&= \nu  \sum_{i=1}^d |v_{h,i}|^2_{H^1(\mathbb{T}_h)} + \nu \eta \sum_{i=1}^d |v_{h,i}|^2_{1} \geq \nu \min \{ 1,\eta \} |{\bm v}_h|_{W_{h}^{CR}}^2. \label{discoe}
\end{align}
\end{rem}

\begin{rem}
In \cite{BarBre14}, the WOPSIP method with $\nu =1$ and $\eta =1$ was proposed for the Stokes equation.
\end{rem}

We can obtain the discrete inf-sup condition with an argument similar to that in \cite[Lemma 7]{Ish24}.

\begin{lem}[Inf-sup stability] \label{lem=infsup}
The nonconforming Stokes element of type $W_{h}^{CR} \times Q_h^{0}$ satisfies the uniform inf-sup stability condition
\begin{align}
\displaystyle
\inf_{q_h \in Q_h^{0}} \sup_{{\bm v}_h \in W_{h}^{CR}} \frac{B_h({\bm v}_h,q_h)}{|{\bm v}_h|_{W_{h}^{CR}} \| q_h \|_{Q_h^{0}}} \geq \beta_* \coloneq \beta. \label{dis=stokes4}
\end{align}
\end{lem}

\subsection{Consistency error estimate}
Let $({\bm u},p) \in W({\bm g}) \times Q$ be the solution to \eqref{stokes2}. Let  ${\bm u}_h \in W_{h,\div}^{CR}$ be the solution to \eqref{dis=stokes3}. For any ${\bm w}_h \in W_{h,\div}^{CR}$, we have that
\begin{align}
\displaystyle
|{\bm u- \bm u_h}|_{W_h^{CR}} &\leq  |{\bm u- \bm w_h}|_{W_h^{CR}}  +  |{\bm w_h-\bm u_h}|_{W_h^{CR}}. \label{dis=stokes5}
\end{align}
We set ${\bm \varphi_h \coloneq \bm w_h- \bm u_h} \in W_{h,\div}^{CR}$. The boundedness and coercivity \eqref{discoe} of $A_h$ and \eqref{dis=stokes3} yield
\begin{align*}
\displaystyle
\nu \min \{ 1,\eta \} |\bm \varphi_h|^2_{W_h^{CR}}
&\leq \nu A_h(\bm \varphi_h,\bm \varphi_h)  \\
&= \nu A_h(\bm w_h - \bm u,\bm \varphi_h) + \nu A_h(\bm u,\bm \varphi_h) - \nu A_h(\bm u_h,\bm \varphi_h) \\
&= \nu A_h(\bm w_h - \bm u, \bm \varphi_h) + \nu A_h(\bm u, \bm \varphi_h) - L_h(\bm \varphi_h) \\
&\leq \nu \max \{1 , \eta \} |\bm w_h - \bm u|_{W_h^{CR}} |\bm \varphi_h|_{W_h^{CR}} + E_h(\bm u) |\bm \varphi_h|_{W_h^{CR}},
\end{align*}
which leads to
\begin{align*}
\displaystyle
|\bm u- \bm u_h|_{W_h^{CR}}
 &\leq \left( 1 + \frac{\max \{1 , \eta \}}{\min \{ 1,\eta \}}\right) \inf_{\bm w_h \in W_{h,\div}^{CR}}  |\bm u- \bm w_h|_{W_h^{CR}} + \frac{1}{\nu \min \{ 1,\eta \}} E_h(\bm u),
\end{align*}
where the consistency term is defined as
\begin{align*}
\displaystyle
 E_h(\bm u) \coloneq \sup_{\bm v_h \in  W_{h,\div}^{CR}} \frac{|\nu A_h(\bm u,\bm v_h) - L_h(\bm v_h)|}{ |\bm v_h|_{W_h^{CR}}}.
\end{align*}
The essential part of the error estimate is the consistency error term.

\color{black}

\begin{lem}[Asymptotic Consistency]  \label{lem=consist}
Assume that $\Omega$ is convex. Let $({\bm u},p) \in (W({\bm g}) \cap H^{2}(\Omega)^d) \times (Q \cap H^{1}(\Omega))$ be the solution to the homogeneous Dirichlet Stokes problem \eqref{stokes2} with ${\bm f} \in L^2(\Omega)^d$ and ${\bm g} \in H^{\frac{1}{2}}(\partial \Omega)^d$ satisfying $\int_{\partial \Omega} {\bm g} \cdot {\bm n} ds = 0$. Let $\{ \mathbb{T}_h\}$ be a family of conformal meshes with the semi-regular property (Assumption \ref{neogeo=assume}). Let $T \in \mathbb{T}_h$ be the element with Conditions \ref{cond1} or \ref{cond2} and satisfy (Type \roman{sone}) in Section \ref{element=cond} when $d=3$. Then,
\begin{align}
\displaystyle
& \frac{1}{\nu \min \{ 1,\eta \}} E_h({\bm u}) \notag \\
 &\quad \leq \frac{c}{\min \{ 1,\eta \}} \left\{ \left( \sum_{i,j=1}^d \sum_{T \in \mathbb{T}_h} h_j^2 \left \| \frac{\partial}{\partial {\bm r}_j} {\bm \nabla} u_i \right \|_{L^2(T)^d}^2 \right)^{\frac{1}{2}} + h \| \varDelta {\bm u} \|_{L^2(\Omega)^d} \right\} \notag \\
&\quad \quad + \frac{c}{\min \{ 1,\eta \}} \left (  h | {\bm u} |_{H^1(\Omega)^d} + h^{ \frac{3}{2}} | {\bm u} |_{H^1(\Omega)^d}^{\frac{1}{2}} \| \varDelta {\bm u} \|_{L^2(\Omega)^d}^{\frac{1}{2}} \right).  \label{dis=stokes6}
\end{align}
\end{lem}

\begin{pf*}
For  $i=1,\ldots,d$, we first have
\begin{align*}
\displaystyle
 \nu \div \mathcal{I}_{h}^{RT} {\bm \nabla} u_i  = \nu \Pi_h^{0} \varDelta u_i.
\end{align*}
From \eqref{disPoi=1} with ${\bm w \coloneq \nu \bm \nabla u_i}$, we have, for any ${\bm v_h} \coloneq (v_{h,1},\ldots,v_{h,d})^{\top} \in W_{h,\div}^{CR}$,
\begin{align}
\displaystyle
&\int_{\Omega} \left( \mathcal{I}_h^{RT} (\nu {\bm \nabla} u_i) \cdot {\bm \nabla_h} v_{h,i} + \nu \Pi_h^{0} \varDelta u_i  v_{h,i} \right) dx = \sum_{F \in \mathcal{F}_h^{\partial}} \int_{F} ( \nu {\bm \nabla} u_i \cdot n_F) \Pi_F^0 v_{h,i} ds.  \label{dis=stokes7}
\end{align}
Thus, using \eqref{dis=stokes7},
\begin{align*}
\displaystyle
&\nu a_h(u_{i} , v_{h,i}) - \ell^i_h(\bm v_h) \\
&\quad = \nu ( \bm \nabla u_i - \mathcal{I}_h^{RT} \bm \nabla u_i , \bm \nabla_h v_{h,i}  )_{\Omega} - \left( - \nu \varDelta u_i + \frac{\partial p}{\partial x_i} , ( \mathcal{I}_{h0}^{RT} \bm v_h )_i \right)_{\Omega} \\
&\quad \quad  - \nu \left( \Pi_h^{0} \varDelta u_i , v_{h,i} -  (\mathcal{I}_{h0}^{RT} \bm v_h )_i \right)_{\Omega} - \nu \left(  \Pi_h^{0} \varDelta u_i ,  (\mathcal{I}_{h0}^{RT} \bm v_h )_i  \right)_{\Omega} \\
&\quad \quad +  \nu \sum_{F \in \mathcal{F}_h^{\partial}} \int_{F} ( \bm \nabla u_i \cdot \bm n_F) \Pi_F^0 v_{h,i} ds \\
&\quad = \nu ( \bm \nabla u_i - \mathcal{I}_h^{RT} \bm \nabla u_i , \bm \nabla_h v_{h,i}  )_{\Omega} + \nu \left( \varDelta u_i  -  \Pi_h^{0} \varDelta u_i, ( \mathcal{I}_{h0}^{RT} \bm v_h )_i \right)_{\Omega} \\
&\quad \quad  - \nu \left( \Pi_h^{0} \varDelta u_i , v_{h,i} -  (\mathcal{I}_{h0}^{RT} \bm v_h )_i \right)_{\Omega}  - \left( \frac{\partial p}{\partial x_i} , ( \mathcal{I}_{h0}^{RT} \bm v_h )_i \right)_{\Omega} \\
&\quad \quad + \nu \sum_{F \in \mathcal{F}_h^{\partial}} \int_{F} ( \bm \nabla u_i \cdot \bm n_F) \Pi_F^0 v_{h,i} ds,
\end{align*}
which leads to
\begin{align*}
\displaystyle
&\nu A_h(\bm u,\bm v_h) - L_h(\bm v_h) \\
&\quad =  \nu \sum_{i=1}^d ( \bm \nabla u_i - \mathcal{I}_h^{RT} \bm \nabla u_i , \bm \nabla_h v_{h,i}  )_{\Omega} + \nu \sum_{i=1}^d \left( \varDelta u_i  -  \Pi_h^{0} \varDelta u_i, ( \mathcal{I}_{h0}^{RT} \bm v_h )_i \right)_{\Omega} \\
&\quad \quad  - \nu \sum_{i=1}^d \left( \Pi_h^{0} \varDelta u_i , v_{h,i} -  (\mathcal{I}_{h0}^{RT} \bm v_h )_i \right)_{\Omega} - \int_{\Omega} ( ( \mathcal{I}_{h0}^{RT} \bm v_h ) \cdot \bm \nabla )p dx \\
&\quad \quad + \nu \sum_{i=1}^d  \sum_{F \in \mathcal{F}_h^{\partial}} \int_{F} ( \bm \nabla u_i \cdot \bm  n_F) \Pi_F^0 v_{h,i} ds \\
&\quad \eqcolon K_1 + K_2 +K_3 +K_4 + K_5.
\end{align*}
\color{black}
Using the H\"older inequality, the Cauchy--Schwarz inequality, and the RT interpolation error \eqref{RT5}, the term $K_1$ is estimated as
\begin{align*}
\displaystyle
|K_1| 
&\leq c \nu  \left( \left( \sum_{T \in \mathbb{T}_h}  \sum_{i,j=1}^d {h}_j^2  \left \| \frac{\partial}{\partial {\bm r}_j} {\bm \nabla} u_i \right \|_{L^2(T)^d}^2 \right)^{\frac{1}{2}} +  h   \|  \varDelta {\bm u} \|_{L^2(\Omega)^d} \right)  | {\bm v}_{h}|_{W_h^{CR}}.
\end{align*}
Using the H\"older inequality, the Cauchy--Schwarz inequality, the stability of $\Pi_h^0$, the estimate \eqref{L2ortho} and \eqref{RT5}, the term $K_2$ is estimated as
\begin{align*}
\displaystyle
|K_2|
&\leq \nu \sum_{i=1}^d \int_{\Omega} |  \varDelta u_i  - \Pi_h^{0} \varDelta u_i |  |( \mathcal{I}_{h0}^{RT} {\bm v}_{h})_i - v_{h,i}| dx \\
&\quad + \nu \sum_{i=1}^d \int_{\Omega} |  \varDelta u_i  - \Pi_h^{0} \varDelta u_i |  | v_{h,i} - \Pi_h^0 v_{h,i}| dx \\
&\leq c \nu h \| \varDelta {\bm u} \|_{L^2(\Omega)^d} | {\bm v}_h|_{W_h^{CR}},
\end{align*}
where we used the fact that
\begin{align*}
\displaystyle
 \int_{\Omega} (  \varDelta u_i  - \Pi_h^{0} \varDelta u_i ) \Pi_h^0 v_{h,i} dx = 0.
\end{align*}
The RT interpolation error \eqref{RT5}, the H\"older's, Cauchy--Schwarz inequalities and the stability of $\Pi_h^0$ yields
\begin{align*}
\displaystyle
|K_3|
&\leq c \nu h \| \varDelta {\bm u} \|_{L^2(\Omega)^d} | {\bm v}_h|_{W_h^{CR}}. 
\end{align*}
Integration by parts,
\begin{align*}
\displaystyle
K_4 =- \int_{\Omega} ( (\mathcal{I}_{h0}^{RT} {\bm v}_{h}) \cdot {\bm \nabla}) p dx = 0.
\end{align*}
where $p \in H_*^1(\Omega)$ and $\mathcal{I}_{h0}^{RT} {\bm v}_{h} \in V_{h0}^{RT}$. Using inequality \eqref{disPoi=2}, the term $K_5$ is estimated as
\begin{align*}
\displaystyle
|K_5|
&\leq c \nu \left(  h | {\bm u} |_{H^1(\Omega)^d} + h^{ \frac{3}{2}} | {\bm u} |_{H^1(\Omega)^d}^{\frac{1}{2}} \| \varDelta {\bm u} \|_{L^2(\Omega)^d}^{\frac{1}{2}} \right)   | {\bm v}_{h} |_{W_h^{CR}}.
\end{align*}
The target estimate \eqref{dis=stokes6} is deduced by gathering the above inequalities.
\qed
\end{pf*}

\begin{rem}[$L^2$-error estimate for the velosity]
The $L^2$-error estimate of a well-balanced scheme on shape-regular meshes was proven in \cite{LinMer17}. However, $L^2$-error analysis of anisotropic meshes remains open.
\end{rem}

\subsection{Error estimate for pressure}
This section provides the $ L^2$-error estimate for the pressure.

\begin{thr}
Assume that $\Omega$ is convex. Let $({\bm u},p) \in (W({\bm g}) \cap H^{2}(\Omega)^d) \times (Q \cap H^{1}(\Omega))$ be the solution of the homogeneous Dirichlet Stokes problem \eqref{stokes2} with ${\bm f} \in L^2(\Omega)^d$ and ${\bm g} \in H^{\frac{1}{2}}(\partial \Omega)^d$ satisfying $\int_{\partial \Omega} {\bm g} \cdot {\bm n} ds = 0$. Let $\{ \mathbb{T}_h\}$ be a family of conformal meshes with the semi-regular property (Assumption \ref{neogeo=assume}). Let $T \in \mathbb{T}_h$ be the element with Conditions \ref{cond1} or \ref{cond2} and satisfy (Type \roman{sone}) in Section \ref{element=cond} when $d=3$. Let $({\bm u}_h,p_h) \in W_h^{CR} \times Q_h^0$ be the solution to \eqref{dis=stokes2}. Then,
\begin{align}
\displaystyle
\| p - p_h \|_{Q_h^0}
&\leq  c \nu \left\{ \left( \sum_{i,j=1}^d \sum_{T \in \mathbb{T}_h} h_j^2 \left \| \frac{\partial}{\partial {\bm r}_j} {\bm \nabla} u_i \right \|_{L^2(T)^d}^2 \right)^{\frac{1}{2}} + h \| \varDelta {\bm u} \|_{L^2(\Omega)^d} \right\} \notag \\
&\quad \quad + c \nu \left (  h | {\bm u} |_{H^1(\Omega)^d} + h^{ \frac{3}{2}} | {\bm u} |_{H^1(\Omega)^d}^{\frac{1}{2}} \| \varDelta {\bm u} \|_{L^2(\Omega)^d}^{\frac{1}{2}} \right) \notag \\
&\quad \quad +c \left( \sum_{j=1}^d \sum_{T \in \mathbb{T}_h} h_j^2 \left\| \frac{\partial p}{\partial {\bm r}_j} \right\|_{L^{2}(T)}^2 \right)^{\frac{1}{2}}.  \label{dis=stokes8}
\end{align}
\end{thr}

\begin{pf*}
To estimate the pressure error $ \| p - p_h \|_{Q_h^0}$, we use the inf-sup stability relation \eqref{dis=stokes4}.  For an arbitrary $q_h \in Q_h^0$, it follows that
\begin{align}
\displaystyle
\| p - p_h \|_{Q_h^0}
&\leq  \| p - q_h \|_{Q_h^0} + \| q_h - p_h \|_{Q_h^0} \notag \\
&\leq  \| p - q_h \|_{Q_h^0} + \frac{1}{\beta_*}  \sup_{{\bm v}_h \in W_{h}^{CR}} \frac{B_h({\bm v}_h, q_h - p_h)}{|{\bm v}_h|_{W_{h}^{CR}}} \notag \\
&\leq  \| p - q_h \|_{Q_h^0} \notag \\
&\quad + \frac{1}{\beta_*}  \sup_{{\bm v}_h \in W_{h}^{CR}} \frac{B_h({\bm v}_h, q_h - p)}{|{\bm v}_h|_{W_{h}^{CR}}} + \frac{1}{\beta_*}  \sup_{{\bm v}_h \in W_{h}^{CR}} \frac{B_h({\bm v}_h, p - p_h)}{|{\bm v}_h|_{W_{h}^{CR}}} \notag \\
&\leq \left( 1 + \frac{1}{\beta_*} \right)  \| p - q_h \|_{Q_h^0} + \frac{1}{\beta_*}  \sup_{{\bm v}_h \in W_{h}^{CR}} \frac{B_h({\bm v}_h, p - p_h)}{|{\bm v}_h|_{W_{h}^{CR}}}. \label{dis=stokes9}
\end{align}
From ${\bm f = - \nu \varDelta \bm u + \bm \nabla p}$, we obtain
\begin{align}
\displaystyle
- L_h({\bm v}_h)
&= \nu \int_{\Omega} \varDelta {\bm u} \cdot \mathcal{I}_{h0}^{RT} {\bm v}_h dx - \int_{\Omega} {\bm \nabla} p \cdot \mathcal{I}_{h0}^{RT} {\bm v}_h dx \notag \\
&\quad - \nu \eta \sum_{F \in \mathcal{F}_h^{\partial}} \kappa_{F(1)} \int_F \Pi_F^0 {\bm g} \cdot \Pi_F^0 {\bm v}_h ds. \label{dis=stokes10}
\end{align}
Using the relation between the RT and CR finite elements (see Lemma \ref{Nit=lem1})  yields
\begin{align}
\displaystyle
\nu A_h(\bm u, \bm v_h)  &= \nu \sum_{i=1}^d \sum_{T \in \mathbb{T}_h} \int_{T}( \bm  \nabla u_{i} - \mathcal{I}_{h}^{RT} \bm \nabla u_i) \cdot \nabla v_{h,i} dx - \nu \int_{\Omega} \Pi_h^0 \varDelta \bm u \cdot \bm v_h dx \notag \\
&\quad +  \nu \eta \sum_{F \in \mathcal{F}_h^{\partial}} \kappa_{F(1)} \int_F \Pi_F^0 \bm u \cdot \Pi_F^0 \bm v_h ds \notag \\
&\quad +  \nu \sum_{i=1}^d \sum_{F \in \mathcal{F}_h^{\partial}} \int_{F} ( \bm \nabla u_i \cdot \bm n_F) \Pi_F^0 v_{h,i} ds. \label{dis=stokes11}
\end{align}
Using the Gauss--Green formula, we obtain
\begin{align}
\displaystyle
- \int_{\Omega} \bm \nabla p \cdot \mathcal{I}_{h0}^{RT} \bm v_h dx
&= \int_{\Omega} p \Pi_h^0 \divh v_h dx = - B_h(\bm v_h,p). \label{dis=stokes12}
\end{align}
\color{black}
Here, we used $\Pi_h^0 \divh {\bm v}_h =  \divh {\bm v}_h$ because $\div {\bm v}_h$ is constant on $T$. From \eqref{dis=stokes2a}, \eqref{dis=stokes10}, \eqref{dis=stokes11} and  \eqref{dis=stokes12}, we have
\begin{align*}
\displaystyle
B_h(\bm v_h, p - p_h)
&= B_h(\bm v_h, p ) + \nu A_h(\bm u_h - \bm u,\bm v_h) + \nu A_h(\bm u,\bm v_h)  - L_h(\bm v_h) \\
&= \nu A_h(\bm u_h - \bm u, \bm v_h) + \nu \sum_{i=1}^d \sum_{T \in \mathbb{T}_h} \int_{T}( \bm \nabla u_{i} - \mathcal{I}_{h}^{RT} \bm \nabla u_i) \cdot \bm \nabla v_{h,i} dx \\
&\quad + \nu \int_{\Omega} \varDelta \bm u \cdot \mathcal{I}_{h0}^{RT} \bm v_h dx - \nu \int_{\Omega} \Pi_h^0 \varDelta \bm u \cdot \bm v_h dx \\
&\quad +  \nu \sum_{i=1}^d \sum_{F \in \mathcal{F}_h^{\partial}} \int_{F} ( \bm \nabla u_i \cdot \bm n_F) \Pi_F^0 v_{h,i} ds \\
&\eqcolon L_1 + L_2 + L_3 + L_4 + L_5.
\end{align*}
\color{black}
Using the CR interpolation error \eqref{CRint} and Lemma \ref{lem=consist}, the term  $L_1$ is estimated as follows:
\begin{align*}
\displaystyle
|L_1| &\leq c \nu |{\bm u - \bm u_h}|_{W_h^{CR}} |{\bm v_h}|_{W_h^{CR}} \\
&\leq  c \nu \left\{ \left( \sum_{i,j=1}^d \sum_{T \in \mathbb{T}_h} h_j^2 \left \| \frac{\partial}{\partial {\bm r}_j} {\bm \nabla} u_i \right \|_{L^2(T)^d}^2 \right)^{\frac{1}{2}} + h \| \varDelta {\bm u} \|_{L^2(\Omega)^d} \right\}  |{\bm v}_h|_{W_h^{CR}} \notag \\
&\quad \quad + c \nu \left (  h | {\bm u} |_{H^1(\Omega)^d} + h^{ \frac{3}{2}} | {\bm u} |_{H^1(\Omega)^d}^{\frac{1}{2}} \| \varDelta {\bm u} \|_{L^2(\Omega)^d}^{\frac{1}{2}} \right)  |{\bm v}_h|_{W_h^{CR}}.
\end{align*}
As proof of Lemma \ref{lem=consist}, $L_2+L_3+L_4+L_5$ is estimated as
\begin{align*}
\displaystyle
&|L_2+L_3+L_4+L_5| \\
&\quad \leq c \nu \left\{ \left( \sum_{i,j=1}^d \sum_{T \in \mathbb{T}_h} h_j^2 \left \| \frac{\partial}{\partial {\bm r}_j} {\bm \nabla} u_i \right \|_{L^2(T)^d}^2 \right)^{\frac{1}{2}} + h \| \varDelta {\bm u} \|_{L^2(\Omega)^d} \right\}  |{\bm v}_h|_{W_h^{CR}} \notag \\
&\quad \quad + c \nu  \left (  h | {\bm u} |_{H^1(\Omega)^d} + h^{ \frac{3}{2}} | {\bm u} |_{H^1(\Omega)^d}^{\frac{1}{2}} \| \varDelta {\bm u} \|_{L^2(\Omega)^d}^{\frac{1}{2}} \right)  |{\bm v}_h|_{W_h^{CR}}.
\end{align*}
From \eqref{L2ortho}, it holds that
\begin{align*}
\displaystyle
\inf_{q_h \in Q_h^0} \| p - q_h \|_{Q_h^0} 
&\leq \| p - \Pi_{h}^0 p \|_{Q_h^0} \leq c \left( \sum_{j=1}^d \sum_{T \in \mathbb{T}_h} h_j^2 \left\| \frac{\partial p}{\partial {\bm r}_j} \right\|_{L^{2}(T)}^2 \right)^{\frac{1}{2}}. 
\end{align*}
We gather the above inequalities and obtain the target estimate \eqref{dis=stokes8}.
\qed
\end{pf*}

\section{Numerical experiments} \label{sec:8}
\setcounter{section}{7} \setcounter{equation}{0} 
Let $d=2$ and $\Omega := (0,1)^2$. The function ${\bm f}$ of the Stokes equation \eqref{stokes1},
\begin{align}
\displaystyle
- \nu \varDelta {\bm u} + {\bm \nabla} p = {\bm f} \quad \text{in $\Omega$}, \quad \div {\bm u}  = 0 \quad \text{in $\Omega$}, \quad {\bm u = \bm g} \quad \text{on $\partial \Omega$}, \label{num1}
\end{align}
is given as satisfying exact solutions $(u,p)$, where ${\bm g}:\Omega \to \mathbb{R}^2$ is a given function and ${\bm g} \in H^{\frac{1}{2}}(\partial \Omega)^2$ satisfies $\int_{\partial \Omega} {\bm g \cdot \bm n} ds = 0$. 

If an exact solution $u$ is known, then the errors ${\bm e_N := \bm u - \bm u_N}$ and ${\bm e_{2N} := \bm u - \bm u_{2N}}$ are computed numerically for the {two division numbers $N$ and $2N$.} The convergence indicator $r$ is defined as follows:
\begin{align*}
\displaystyle
r = \frac{1}{\log(2)} \log \left( \frac{\| {\bm e}_N \|_X}{\| {\bm e}_{2N} \|_X} \right).
\end{align*}
We compute the convergence order concerning the norms defined by
\begin{align*}
\displaystyle
Err(W_h^{CR}) &:= \frac{| {\bm u - \bm u_N} |_{W_h^{CR}}}{|{\bm u}|_{W_h^{CR}}}, \quad Err(L^2):= \frac{\| {\bm u - \bm u_N} \|_{L^2(\Omega)^d}}{\| {\bm u} \|_{L^2(\Omega)^d}},\\
Err(Q_h^0) &:= \frac{\| p - p_N \|_{Q_h^0}}{\|p\|_{Q_h^0}}.
\end{align*}

\begin{figure}[htbp]
   \includegraphics[keepaspectratio, scale=0.15]{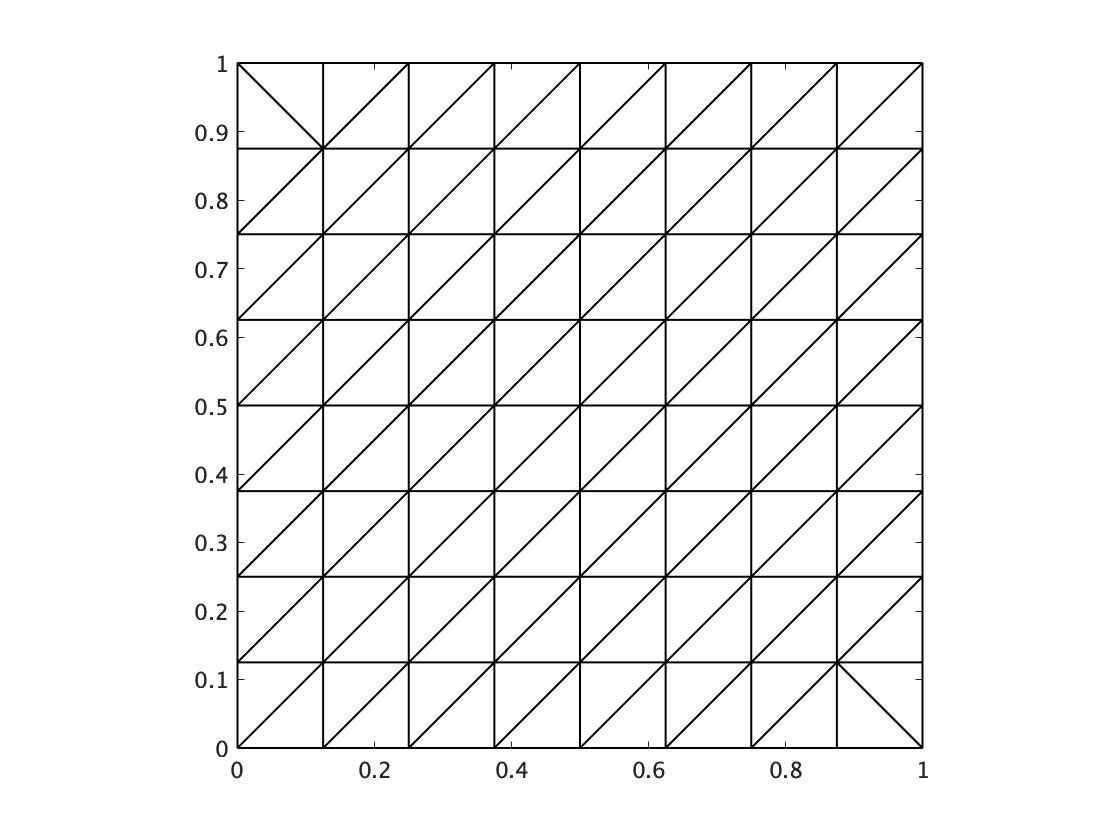}
    \caption{(Mesh \Roman{lone})  $N=8$}
     \label{fig1}
\end{figure}

\begin{figure}[htbp]
  \begin{minipage}[b]{0.45\linewidth}
    \centering
    \includegraphics[keepaspectratio, scale=0.15]{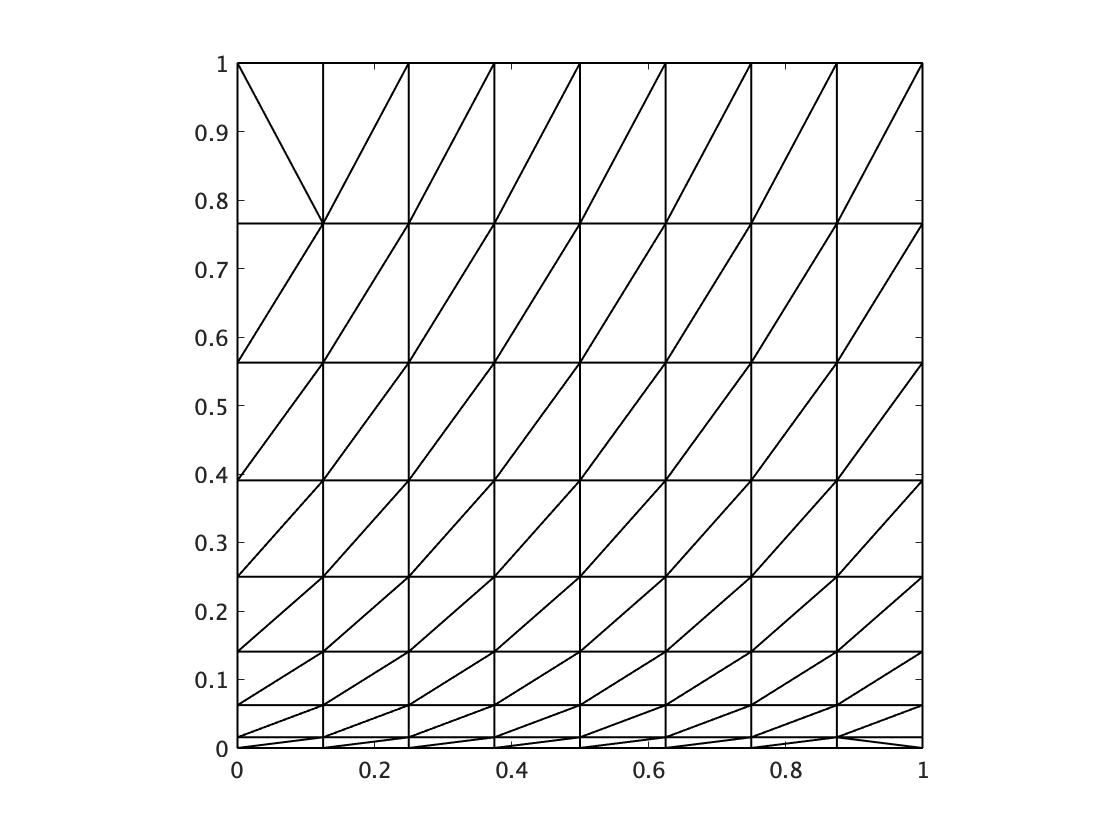}
    \caption{(Mesh \Roman{ltwo}) $N=8$}
     \label{fig00}
  \end{minipage}
  \begin{minipage}[b]{0.45\linewidth}
    \centering
    \includegraphics[keepaspectratio, scale=0.15]{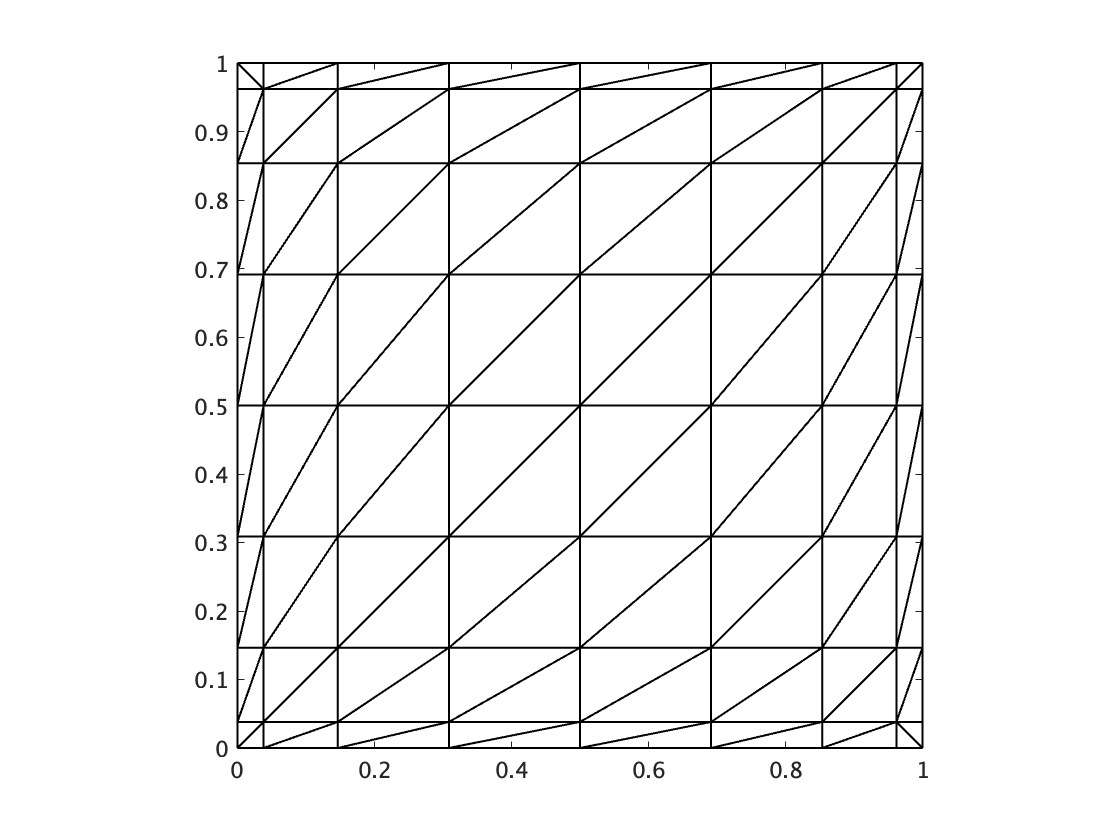}
    \caption{(Mesh \Roman{lthree})  $N=8$}
     \label{fig002}
  \end{minipage}
\end{figure}

Let $N \in \{ 16,32,64,128,256,512\}$ be the division number of each side of the bottom and the height edges of $\Omega$. The following mesh partitions are considered.  Let $(x_1^i, x_2^i)^{\top}$ and $i \in \mathbb{N}$ be the grip points of the triangulations $\mathbb{T}_h$ defined as follows: 
\begin{description}
 \item[(Mesh \Roman{lone}) Standard mesh (Fig. \ref{fig1})] 
\begin{align*}
\displaystyle
x_1^i \coloneq \frac{i}{N}, \quad x_2^i \coloneq \frac{i}{N}, \quad  i \in \{0, \ldots , N \}.
\end{align*}
  \item[(Mesh \Roman{ltwo})]  Graded mesh (Fig. \ref{fig00})
\begin{align*}
\displaystyle
x_1^i \coloneq \frac{i}{N}, \quad x_2^i \coloneq \left ( \frac{i}{N} \right)^{2}, \quad  i \in \{0, \ldots, N \}.
\end{align*}
 \item[(Mesh \Roman{lthree})  Anisotropic mesh which comes from \cite{CheLiuQia10} (Fig. \ref{fig002})]
\begin{align*}
\displaystyle
x_1^i \coloneq \frac{1}{2}\left( 1 - \cos \left( \frac{i \pi}{N} \right) \right), \quad x_2^i \coloneq\frac{1}{2}\left( 1 - \cos \left( \frac{i \pi}{N} \right) \right), \quad  i \in \{0, \ldots, N \}.
\end{align*}
\end{description}
(Mesh \Roman{lone}) and (Mesh \Roman{ltwo}) were used  Eexample 1 in Section \ref{numerical1}. (Mesh \Roman{lone}) and (Mesh \Roman{ltwo}) were used Example 2 in Section \ref{numerical2}. In \cite[Tables 1 and 6]{Ish24d}, (Mesh \Roman{ltwo}) and (Mesh \Roman{lthree}) do not satisfy the shape-regular condition but satisfy the semiregularity mesh condition defined in \eqref{NewGeo}, which is equivalent to the maximum angle condition. 

For the computation, we used the CG method without preconditioners and quadrature of the five orders \cite[p. 85, Table 30.1]{ErnGue21b}.

\subsection{Example 1} \label{numerical1}
The first example shows the results of numerical tests on anisotropic meshes. The function ${\bm f}$ of the Stokes equation \eqref{num1} with $\nu = 1.0, 10^{-5}$ is given such that the exact solution is
\begin{align*}
\displaystyle
\begin{pmatrix}
 u_1  \\
  u_2
\end{pmatrix}
&= 
\begin{pmatrix}
 \sin \pi x_1 \cos \pi x_2 \\
  - \cos \i x_1 \sin \pi x_2
\end{pmatrix},\quad
p =  \sin \pi x_1 \cos \pi x_2.
\end{align*}
The numerical results of the scheme \eqref{dis=stokes2} for Example 1 are listed in Tables \ref{table4} -- \ref{table4dd}. In Example 1, the following numerical results were obtained.

\begin{description}
  \item[(R1)] Tables \ref{table4}, \ref{table4b} and \ref{table4bb} are the results under the shape-regularity condition. It is observed that the optimal order can be obtained by tuning the appropriate parameter $\eta$ from Tables \ref{table4} and \ref{table4b}. It is observed that the error increases if parameter $\eta$ is not selected appropriately (Table \ref{table4bb}).
  \item[(R2)] Tables \ref{table4c}, \ref{table4d} and \ref{table4dd} are the results under the semi-regularity condition. The same trends as (R1) in Tables \ref{table4}, \ref{table4b} and \ref{table4bb} were generally observed. It was confirmed that optimal order convergence can be obtained even if the shape-regularity mesh condition is relaxed.
\end{description}
Numerical calculations for the scheme \eqref{dis=stokes2} were performed to confirm that the theoretical results were satisfactory. 

\begin{table}[h]
\caption{(Mesh \Roman{lone}), Example 1 with $\nu = 1$, $\eta = 1$}
\centering
\begin{tabular}{l|l|l|l|l|l|l|l} \hline
$N$ &  $h $ & $Err(W_h^{CR})$ & $r$ & $Err(L^2)$  & $r$ & $Err(Q_h)$ & $r$ \\ \hline \hline
16 & 8.84e-02  &   1.20684e-01   &  &9.93092e-03   &   & 1.01648e-01 &   \\
32 & 4.42e-02  & 5.99648e-02    & 1.01 & 2.15748e-03   & 2.02  & 3.27172e-02 & 1.64 \\
64 &  2.21e-02   & 2.96676e-02  & 1.02 & 5.18816e-04 &  2.06 & 1.63615e-02  & 1.00 \\
128 & 1.10e-02    &  9.48850e-03   &  1.64  &7.93848e-05  &    2.71  & 8.18111e-03 & 1.00 \\
256 &   5.52e-03    &  4.71422e-03    & 1.01   &   1.98438e-05  &  2.00   & 4.09060e-03   & 1.00  \\
512 &  2.76e-03     &  2.35001e-03    &  1.00  & 4.96081e-06    & 2.00    & 2.04531e-03   & 1.00 \\
\hline
\end{tabular}
\label{table4}
\end{table}

\begin{table}[h]
\caption{(Mesh \Roman{lone}), Example 1 with $\nu = 10^{-5}$, $\eta = 10^5$}
\centering
\begin{tabular}{l|l|l|l|l|l|l|l} \hline
$N$ &  $h $ & $Err(W_h^{CR})$ & $r$ & $Err(L^2)$  & $r$ & $Err(Q_h)$ & $r$ \\ \hline \hline
16 & 8.84e-02 &   6.79729e-02  &  &4.55688e-03  &   &  6.53786e-02  &   \\
32 &  4.42e-02   &  3.75622e-02  & 0.86  &  1.26946e-03 & 1.84   &  3.27170e-02  & 1.00 \\
64 &    2.21e-02   &1.87647e-02   & 1.00 & 3.17457e-04  & 2.00   &1.63615e-02  & 1.00 \\
128 & 1.10e-02   & 9.37774e-03   & 1.00  & 7.93702e-05 &  2.00 & 8.18111e-03  & 1.00 \\
256 &  5.52e-03    &  4.68765e-03    & 1.00   &   1.98429e-05  &   2.00 & 4.09060e-03   & 1.00 \\
512 &   2.76e-03    &  2.34351e-03    &  1.00  & 4.96076e-06    &   2.00  &  2.04531e-03  & 1.00 \\
\hline
\end{tabular}
\label{table4b}
\end{table}

\begin{table}[h]
\caption{(Mesh \Roman{lone}), Example 1 with $\nu = 10^{-5}$, $\eta = 1$}
\centering
\begin{tabular}{l|l|l|l|l|l|l|l} \hline
$N$ &  $h $ & $Err(W_h^{CR})$ & $r$ & $Err(L^2)$  & $r$ & $Err(Q_h)$ & $r$ \\ \hline \hline
16 &8.84e-02&  7.01415e+02  &  & 4.79195e+01   &   &  6.53664e-02&   \\
32 &  4.42e-02 &  2.48543e+02  & 1.50 & 6.04345e+00  &  2.99 & 3.27125e-02 & 1.00 \\
64 &    2.21e-02   & 8.79111e+01  & 1.50 & 7.56851e-01  &  3.00 & 1.63608e-02  & 1.00 \\
128 &   1.10e-02     &  3.10841e+01  & 1.50 & 9.38478e-02 &  3.01 &  8.18103e-03  & 1.00 \\
256 &   5.52e-03   & 1.09901e+01      & 1.50   &   1.11251e-02  &  3.08   &  4.09075e-03  & 1.00 \\
512 &  2.76e-03    &   3.88599    &  1.50  & 1.33775e-03   &  3.06   & 2.06576e-03  &  0.99\\
\hline
\end{tabular}
\label{table4bb}
\end{table}

\begin{table}[h]
\caption{(Mesh \Roman{ltwo}), Example 1 with $\nu = 1$, $\eta = 1$}
\centering
\begin{tabular}{l|l|l|l|l|l|l|l} \hline
$N$ &  $h $ & $Err(W_h^{CR})$ & $r$ & $Err(L^2)$  & $r$ & $Err(Q_h)$ & $r$ \\ \hline \hline
16 & 1.36e-01 &    1.50239e-01 &  & 1.58750e-02  &   & 1.35046e-01  &   \\
32 &  6.90e-02 &   7.41164e-02  & 1.02  & 3.36663e-03  &  2.24 & 4.00647e-02  & 1.75 \\
64 &   3.47e-02   & 2.45175e-02  & 1.60 & 5.01453e-04  &  2.75 &  2.00380e-02& 1.00 \\
128 & 1.74e-02  &  1.14880e-02   &  1.10 &   1.24652e-04  & 2.01   &  1.00197e-02 & 1.00 \\
256 &  8.72e-03      &    5.74067e-03   & 1.00   & 3.11643e-05   & 2.00    & 5.00993e-03  &1.00  \\
512 &    4.36e-03    &  2.86994e-03    &  1.00  & 7.79116e-06   &   2.00  & 2.50498e-03   &1.00  \\
\hline
\end{tabular}
\label{table4c}
\end{table}

\begin{table}[h]
\caption{(Mesh \Roman{ltwo}), Example 1 with  $\nu = 10^{-5}$, $\eta = 10^5$}
\centering
\begin{tabular}{l|l|l|l|l|l|l|l} \hline
$N$ &  $h $ & $Err(W_h^{CR})$ & $r$ & $Err(L^2)$  & $r$ & $Err(Q_h)$ & $r$ \\ \hline \hline
16 & 1.36e-01&   8.46574e-02  &  & 7.34737e-03  &   &  8.00174e-02  &   \\
32 &  6.90e-02 & 4.60690e-02    & 0.88 &  1.99401e-03  &  1.88 & 4.00645e-02  & 1.00 \\
64 &   3.47e-02   & 2.29559e-02 & 1.00 &  4.98496e-04  &  2.00 &  2.00380e-02  & 1.00 \\
128 & 1.74e-02  & 1.14789e-02    & 1.00  & 1.24650e-04 &  2.00 &  1.00197e-02  & 1.00 \\
256 &   8.72e-03   &  5.73959e-03    &  1.00  & 3.11642e-05  &  2.00   & 5.00993e-03  & 1.00  \\
512 &   4.36e-03   &  2.86981e-03     &  1.00  & 7.79116e-06   &   2.00  &   2.50498e-03  & 1.00 \\
\hline
\end{tabular}
\label{table4d}
\end{table}

\begin{table}[h]
\caption{(Mesh \Roman{ltwo}), Example 1 with  $\nu = 10^{-5}$, $\eta = 1$}
\centering
\begin{tabular}{l|l|l|l|l|l|l|l} \hline
$N$ &  $h $ & $Err(W_h^{CR})$ & $r$ & $Err(L^2)$  & $r$ & $Err(Q_h)$ & $r$ \\ \hline \hline
16 & 1.36e-01&  8.62497e+02   &  & 8.72047e+01 &   &  6.04737e-01  &   \\
32 &  6.90e-02 & 3.05103e+02   & 1.50 & 1.14315e+01  &  2.93 & 6.18889e-01  & -0.03 \\
64 &   3.47e-02   & 1.07618e+02   & 1.50 &  1.45683e+00  &  2.97 &  6.27180e-01  & -0.02 \\
128 & 1.74e-02  &  3.79818e+01   & 1.50 & 1.76028e-01  &  3.05 &  6.31746e-01 & -0.01 \\
256 &   8.72e-03   &  1.34116e+01    & 1.50   & 2.20423e-02   &  3.00   &  6.31677e-01 & 1.57e-04  \\
512 &   4.36e-03   &  6.13881    &  1.13  & 1.89780e-04   &  6.86   & 2.50498e-03   & 7.98  \\
\hline
\end{tabular}
\label{table4dd}
\end{table}

\subsection{Example 2} \label{numerical2}
The second example demonstrates the robustness of the proposed scheme against large irrotational body forces on anisotropic mesh partitions. The exact solution $(u,p)$ to \eqref{num1} with velocity and pressure is given by
\begin{align*}
\displaystyle
\begin{pmatrix}
 u_1  \\
  u_2
\end{pmatrix}
&=
\begin{pmatrix}
 - \left( x_2 - \frac{1}{2} \right)  \\
x_1 - \frac{1}{2}
\end{pmatrix}, \\
p &= 10^5 (1 - x_2)^3 - \frac{10^5}{4}.
\end{align*}
Notably, $\int_{\Omega} p dx = 0$. By setting $\nu = 1$, the force is exactly irrotational, that is, 
\begin{align*}
\displaystyle
\begin{pmatrix}
 f_1  \\
 f_2
\end{pmatrix}
&=
\begin{pmatrix}
0 \\
- 3 \times 10^5 (1 - x_2)^2
\end{pmatrix}.
\end{align*}
The numerical results of the scheme \eqref{dis=stokes2} for Example 2 are listed in Tables \ref{table7a} and \ref{table7}. In Example 2, the following numerical results were obtained.

\begin{description}
  \item[(R3)] Table \ref{table7a} is the results under the shape-regularity condition. The optimal order can be obtained by tuning the appropriate parameter $\eta$. 
  \item[(R4)]  Table \ref{table7} lists the results for the semi-regularity condition. The optimal order can also be obtained by tuning the appropriate parameter $\eta$.
  \item[(R5)]  By observing  Table \ref{table7}, the error is slightly reduced when using the anisotropic mesh in Example 2.
\end{description}
Numerical calculations for the scheme \eqref{dis=stokes2} were performed to confirm that the theoretical results were satisfactory. 

\begin{table}[h]
\caption{(Mesh \Roman{lone}), Example 2 with  $\nu = 1$, $\eta = 10^5$}
\centering
\begin{tabular}{l|l|l|l|l|l|l|l} \hline
$N$ &  $h $ & $Err(W_h^{CR})$ & $r$ & $Err(L^2)$  & $r$ & $Err(Q_h)$ & $r$ \\ \hline \hline
16 & 8.84e-02 &  1.39373e-02    &  & 1.08045e-04   &   & 6.97038e-02  &   \\
32 &  4.42e-02  &  4.90448e-03  & 1.51 & 9.56041e-06 & 3.50   & 3.48586e-02  & 1.00 \\
64 &    2.21e-02   &1.73194e-03   & 1.50 & 8.45263e-07 & 3.50  & 1.74301e-02  & 1.00 \\
128 & 1.10e-02  & 6.12153e-04   & 1.50 &7.47166e-08  & 3.50  &  8.71517e-03 & 1.00 \\
256 &   5.52e-03    & 2.16413e-04     &  1.50  &  6.60420e-09   &   3.50  &  4.35760e-03 & 1.00  \\
512 &   2.76e-03    &   7.65122e-05   &  1.50  &  5.83736e-10   &   3.50  & 2.17880e-03 & 1.00 \\ 
\hline
\end{tabular}
\label{table7a}
\end{table}

\begin{table}[h]
\caption{(Mesh \Roman{lthree}), Example 2 with  $\nu = 1$, $\eta = 10^5$}
\centering
\begin{tabular}{l|l|l|l|l|l|l|l} \hline
$N$ &  $h $ & $Err(W_h^{CR})$ & $r$ & $Err(L^2)$  & $r$ & $Err(Q_h)$ & $r$ \\ \hline \hline
16 &  1.38e-01  & 8.62679e-03   &  & 1.58828e-05  &   & 7.54176e-02  &  \\
32 & 6.93e-02   & 2.11516e-03 & 2.03  & 4.96777e-07   &  5.00 & 3.77681e-02    & 1.00 \\
64 & 3.47e-02    &  5.07297e-04  & 2.06 &  1.49737e-08 & 5.05  &    1.88913e-02 & 1.00 \\
128 &  1.74e-02  &  1.26770e-04 & 2.00 & 4.68359e-10  &  5.00 &    9.44657e-03 & 1.00 \\
256 &  8.68e-03     &   3.16875e-05  &  2.00  &   1.46396e-11  &   5.00  & 4.72340e-03  & 1.00  \\
512 &  4.34e-03     & 7.93343e-06     &  2.00  & 4.57513e-13  &   5.00  &  2.36171e-03 & 1.00 \\
\hline
\end{tabular}
\label{table7}
\end{table}

\section{Conclusions} \label{Conclusions}
\setcounter{section}{9} \setcounter{equation}{0}
In this study, we propose the Nitsche method based on the WOPSIP method for solving the Poisson and Stokes equations. We derive optimal order error estimates using appropriate norms on anisotropic meshes. Our numerical experiments concentrate on the pressure-robust scheme for the Stokes problem with nonhomogeneous Dirichlet boundary conditions. According to \cite{BarBre14}, the WOPSIP method "does not require the tuning of any penalty parameter and is very simple to program." As noted in Remark \ref{dis=coe=rem}, the discrete coercivity remains unaffected by the tuning parameter. As indicated in Tables \ref{table4} and \ref{table4c}, when the diffusion coefficient is not particularly small, the numerical error is not expected to be substantial if the parameter $\eta=1$. Conversely, our numerical experiments suggest that errors can be mitigated by adjusting the penalty in accordance with the problem's characteristics. In \cite[Section 6.3]{Ish24c}, it was observed that numerical errors can be reduced by modifying the tuning parameters. The theoretical understanding of this phenomenon remains unresolved, as the precise method for parameter adjustment is still open. It is an intriguing area of research to extend the error analysis to diffusion-convection-reaction, interface, and Navier-Stokes problems. We plan to investigate these topics in future research.

\color{black}
\section*{Funding}
No funding was received to assist with the preparation of this manuscript.


\appendix
\section*{Additional numerical tests in Section \ref{sec:8}}
In this section, we adopt the following scheme. We aim to find $({\bm u}_h,p_h) \in W_{h}^{CR} \times Q_h^{0}$ such that 
\begin{align*}
\displaystyle
\nu A_h({\bm u_h, \bm v_h}) + B_h({\bm v}_h , p_h) &= L_h^*({\bm v}_h) \quad \forall {\bm v}_h \in W_{h}^{CR}, \tag{A1} \\
B_h({\bm u}_h , q_h) &= 0 \quad \forall q_h \in Q_h^{0}, \tag{A2}
\end{align*}
where $A_h: (H^1(\Omega)^d + W_{h}^{CR}) \times (H^1(\Omega)^d + W_{h}^{CR}) \to \mathbb{R}$ and $B_h: (H^1(\Omega)^d+W_{h}^{CR}) \times Q_h^{0} \to \mathbb{R}$ are defined in Section \ref{new=stokes}, and $L_h*: H^1(\Omega)^d + W_h^{CR} \to \mathbb{R}$ is defined as
\begin{align*}
\displaystyle
L_h^*({\bm v}_h) &\coloneq \sum_{i=1}^d \left\{ ( f_i ,  v_{h,i} )_{\Omega} + \nu \eta  \sum_{F \in \mathcal{F}_h^{\partial}} \kappa_{F(1)} \langle \Pi_F^0  g_{i}, \Pi_F^0 v_{h,i} \rangle_F \right\}.
\end{align*}

\subsection*{Example 1 in Section \ref{numerical1}}
The numerical results of the scheme for Example 1 are listed in Tables 9-12.

\begin{table}[h]
\caption{(Mesh \Roman{lone}), Example 1 with  $\nu = 1$, $\eta = 1$}
\centering
\begin{tabular}{l|l|l|l|l|l|l|l} \hline
$N$ &  $h $ & $Err(W_h^{CR})$ & $r$ & $Err(L^2)$  & $r$ & $Err(Q_h)$ & $r$ \\ \hline \hline
16 & 8.84e-02  &   7.89938e-02  &  &  5.40416e-03 &   &  9.90890e-02 &   \\
32 &   4.42e-02  &  3.87989e-02   & 1.03 &1.25617e-03  & 2.11  & 3.27170e-02  & 1.60 \\
64 &    2.21e-02   & 1.88131e-02   & 1.04 &  3.08832e-04  &  2.02 & 1.63615e-02  & 1.00 \\
128 &   1.10e-02    &  9.47902e-03   & 0.99 & 7.93657e-05  &  1.96 &  8.18111e-03  & 1.00 \\
256 &    5.52e-03   &   4.71303e-03  &  1.01  &  1.98426e-05 &  2.00   &  4.09060e-03 &1.00  \\
512 &   2.76e-03    & 2.34987e-03     & 1.00   & 4.96074e-06   &    2.00 &2.04531e-03   &1.00  \\
\hline
\end{tabular}
\label{table000a}
\end{table}

\begin{table}[h]
\caption{(Mesh \Roman{lone}), Example 1 with  $\nu = 10^{-5}$, $\eta = 10^5$}
\centering
\begin{tabular}{l|l|l|l|l|l|l|l} \hline
$N$ &  $h $ & $Err(W_h^{CR})$ & $r$ & $Err(L^2)$  & $r$ & $Err(Q_h)$ & $r$ \\ \hline \hline
16 & 8.84e-02  & 1.23372e+03   &  & 1.57316e+02  &   &  6.83586e-02 &   \\
32 &   4.42e-02  &  6.30014e+02  & 0.97 &  3.95060e+01 & 2.00  &  3.35578e-02 & 1.03 \\
64 &    2.21e-02   &  3.16058e+02  & 1.00 &  9.66287e+00  & 2.03  & 3.03512e-02 & 0.14 \\
128 &   1.10e-02    &  9.37792e-03   & 15.04 & 7.93702e-05  & 16.89  &  8.18111e-03  & 1.89 \\
256 &   5.52e-03     & 4.68768e-03    & 1.00   & 1.98429e-05   &   2.00  & 4.09060e-03   & 1.00  \\
512 &   2.76e-03    &   2.34352e-03    &  1.00  &   4.96076e-06 &   2.00  &  2.04531e-03   & 1.00  \\
\hline
\end{tabular}
\label{table000b}
\end{table}

\begin{table}[h]
\caption{(Mesh \Roman{ltwo}), Example 1 with  $\nu = 1$, $\eta = 1$}
\centering
\begin{tabular}{l|l|l|l|l|l|l|l} \hline
$N$ &  $h $ & $Err(W_h^{CR})$ & $r$ & $Err(L^2)$  & $r$ & $Err(Q_h)$ & $r$ \\ \hline \hline
16 & 1.36e-01   &   1.00907e-01   &  & 8.98765e-03  &   &  1.23104e-01 &   \\
32 &  6.90e-02 & 4.87473e-02   & 1.05 &  1.93924e-03&  2.21 &  4.00645e-02  & 1.62 \\
64 &   3.47e-02    &  2.35594e-02 & 1.05 & 4.98354e-04  &  1.96 & 2.00380e-02 & 1.00 \\
128 &   1.74e-02   &   1.14868e-02   & 1.04 & 1.24650e-04   & 2.00  & 1.00197e-02 & 1.00 \\
256 &  8.72e-03    &  5.74059e-03    & 1.00    & 3.11642e-05   &   2.00  & 5.00995e-03 & 1.00  \\
512 &  4.36e-03     &   2.86994e-03  &  1.00  &   7.79116e-06 &   2.00  &  2.50498e-03 & 1.00  \\
\hline
\end{tabular}
\label{table000c}
\end{table}

\begin{table}[h]
\caption{(Mesh \Roman{ltwo}), Example 1 with  $\nu = 10^{-5}$, $\eta = 10^5$}
\centering
\begin{tabular}{l|l|l|l|l|l|l|l} \hline
$N$ &  $h $ & $Err(W_h^{CR})$ & $r$ & $Err(L^2)$  & $r$ & $Err(Q_h)$ & $r$ \\ \hline \hline
16 &  1.36e-01  &  1.48472e+03   &  & 2.35808e+02  &   &  6.35185e-01 &   \\
32 &   6.90e-02 & 7.63487e+02  & 0.96  & 5.94611e+01  & 1.99  &  5.04366e-02  & 3.65 \\
64 &    3.47e-02   & 2.29561e-02   & 15.02 & 4.98497e-04 &  16.86 &   2.00380e-02 & 1.33 \\
128 &    1.74e-02     &  1.14789e-02  & 1.00 &  1.24650e-04 &  2.00 &  1.00197e-02  & 1.00 \\
256 &  8.72e-03     & 5.73959e-03     &  1.00  & 3.11642e-05   &   2.00  &  5.00993e-03  & 1.00  \\
512 &     4.36e-03  &  2.86981e-03    &  1.00  &  7.79116e-06 &   2.00  & 2.50498e-03  & 1.00 \\
\hline
\end{tabular}
\label{table000d}
\end{table}

\subsection*{Example 2 in Section \ref{numerical2}}
The numerical results of the scheme for Example 2 are listed in Tables 13, 14.

\begin{table}[h]
\caption{(Mesh \Roman{lone}), Example 2 with  $\nu = 1$, $\eta = 10^5$}
\centering
\begin{tabular}{l|l|l|l|l|l|l|l} \hline
$N$ &  $h $ & $Err(W_h^{CR})$ & $r$ & $Err(L^2)$  & $r$ & $Err(Q_h)$ & $r$ \\ \hline \hline
16 & 8.84e-02 & 1.68337e+03    &  &  1.20814e+02  &   &   6.97035e-02 &   \\
32 &  4.42e-02  &  8.40312e+02  & 1.00 & 2.48685e+01  &  2.28 &  3.48586e-02 & 1.00 \\
64 &    2.21e-02   & 1.74623e-03  & 18.88 & 8.65628e-07 &  24.78 & 1.74301e-02 & 1.00 \\
128 & 1.10e-02  &    6.14707e-04  & 1.51 &  7.56217e-08 & 3.52  &  8.71517e-03 & 1.00 \\
256 &  5.52e-03      &   2.16868e-04   & 1.50   &  6.64431e-09  &  3.51   & 4.35760e-03  & 1.00 \\
512 &    2.76e-03    &  7.65946e-05    &  1.50  & 5.85512e-10   &  3.51   & 2.17880e-03  & 1.00  \\
\hline
\end{tabular}
\label{table007a}
\end{table}

\begin{table}[h]
\caption{(Mesh \Roman{lthree}), Example 2 with  $\nu = 1$, $\eta = 10^5$}
\centering
\begin{tabular}{l|l|l|l|l|l|l|l} \hline
$N$ &  $h $ & $Err(W_h^{CR})$ & $r$ & $Err(L^2)$  & $r$ & $Err(Q_h)$ & $r$ \\ \hline \hline
16 &  1.38e-01  & 1.94675e+03   &  &  1.54829e+02    &   & 7.54191e-02    &  \\
32 & 6.93e-02   &  9.65371e+02   & 1.01  & 3.71343e+01    & 2.06   &    3.77681e-02  & 1.00 \\
64 & 3.47e-02    &  5.16552e-04  & 20.83 &  2.18979e-08   & 30.66  &  1.88913e-02   & 1.00 \\
128 &  1.74e-02  & 1.27937e-04   & 2.01 & 6.78913e-10   &  5.01  & 9.44657e-03    & 1.00 \\
256 &   8.68e-03   & 3.18364e-05     & 2.01   &   2.11249e-11  &  5.01   &  4.72340e-03  & 1.00 \\
512 &   4.34e-03   & 7.94210e-06     &  2.00  & 6.58676e-13  &   5.00  & 2.36171e-03  & 1.00  \\
\hline
\end{tabular}
\label{table007b}
\end{table}

\bigskip

\small
\end{document}